\documentclass[twoside]{article}
\usepackage[accepted]{aistats2019}

\usepackage[round]{natbib}

\usepackage[T1]{fontenc}
\usepackage{amsmath}
\usepackage{amsthm}
\usepackage{amssymb}
\usepackage{algorithm}
\usepackage{algcompatible}
\usepackage{booktabs}
\usepackage{color}
\usepackage{graphicx}
\usepackage{multirow}
\usepackage{enumitem}
\usepackage{url}
\usepackage{subcaption}

\usepackage{todonotes}

\usepackage{anits}

\usepackage{xspace}

\begin{document}

%

%

\twocolumn[

\aistatstitle{Towards Gradient Free and Projection Free Stochastic Optimization}



\aistatsauthor{Anit Kumar Sahu \And Manzil Zaheer \And Soummya Kar }

\aistatsaddress{Carnegie Mellon University\And Google AI \And Carnegie Mellon University} 
]
\vspace{20pt}
\begin{abstract}
This paper focuses on the problem of \emph{constrained} \emph{stochastic} optimization. A zeroth order Frank-Wolfe algorithm is proposed, which in addition to the projection-free nature of the vanilla Frank-Wolfe algorithm makes it gradient free. Under convexity and smoothness assumption, we show that the proposed algorithm converges to the optimal objective function at a rate $O\left(1/T^{1/3}\right)$, where $T$ denotes the iteration count. In particular, the primal sub-optimality gap is shown to have a dimension dependence of $O\left(d^{1/3}\right)$, which is the best known dimension dependence among all zeroth order optimization algorithms with one directional derivative per iteration. For non-convex functions, we obtain the \emph{Frank-Wolfe} gap to be $O\left(d^{1/3}T^{-1/4}\right)$. Experiments on black-box optimization setups demonstrate the efficacy of the proposed algorithm.
\end{abstract}

\section{Introduction}
In this paper, we aim to solve the following stochastic optimization problem:
\begin{equation}
\label{eq:prob-def}
\min_{\mathbf{x}\in\mathcal{C}} f\left(\mathbf{x}\right) = \min_{x\in\mathcal{C}}\mathbb{E}_{\mathbf{y}\sim\mathcal{P}}\left[F\left(\mathbf{x};\mathbf{y}\right)\right],
\end{equation}
where $\mathcal{C}\in\mathbb{R}^{d}$ is a closed convex set.
This problem of stochastic constrained optimization has been a focus of immense interest in the context of convex functions \citep{bubeck2015convex} and non-convex functions especially in the context of deep learning \citep{goodfellow2016deep}. 
Solutions to the problem \eqref{eq:prob-def} can be broadly classified into two classes: algorithms which require a projection at each step, for example, projected gradient descent \citep{bubeck2015convex} and projection free methods such as the Frank-Wolfe algorithm \citep{jaggi2013revisiting}. 
Furthermore, algorithms designed to solve the above optimization problem access various kinds of oracles, i.e., first order oracle (gradient queries) and zeroth order oracle (function queries). 
In this paper, we focus on a stochastic version of projection free method, namely Frank-Wolfe algorithm, with access to a zeroth order oracle.

Derivative free optimization or zeroth order optimization is motivated by settings where the analytical form of the function is not available or when the gradient evaluation is computationally prohibitive. 
Developments in zeroth order optimization has been fueled by various applications ranging from problems in medical science, material science and chemistry~\citep{gray2004optimizing,marsden2008computational,gray2004optimizing,deming1978review,marsden2007trailing}. 
In the context of machine learning, zeroth order methods have been applied to attacks on deep neural networks using black box models~\citep{chen2017zoo}, scalable policy optimization for reinforcement learning~\citep{choromanski2018structured} and optimization with bandit feedback~\citep{bubeck2012regret}.

\begin{table*}[t]
    \centering
    \renewcommand{\arraystretch}{1.2}
    \begin{tabular}{@{}llllr@{}}
    \toprule
         \textbf{Reference} & \textbf{Setting} & \textbf{Memory} & \textbf{Primal Rate} & \textbf{Oracle} \\
    \midrule
         \citet{jaggi2013revisiting}&Det. Convex&-&$O(1/t)$&SFO\\
         \citet{hazan2012projection}&Stoch. Convex&$O(t)$&$O(1/t^{1/2})$&SFO\\
         \citet{mokhtari2018conditional}&Stoch. Convex&$O(1)$&$O(1/t^{1/3})$&SFO\\
         \citet{lacoste2016convergence}&Det. Non-convex&-&$O(1/t^{1/2})$&SFO\\
         \citet{reddi2016stochastic}&Stoch. Non-convex&$O(\sqrt{t})$&$O(1/t^{1/4})$&SFO\\
         RDSA [Theorem \ref{th:2}(1)]&Stoch. Convex&$1$&$O(d^{1/3}/t^{1/3})$&SZO\\
         I-RDSA [Theorem \ref{th:2}(2)]&Stoch. Convex&$m$&$O((d/m)^{1/3}/t^{1/3})$&SZO\\
         KWSA [Theorem \ref{th:2}(3)]&Stoch. Convex&$d$&$O(1/t^{1/3})$&SZO\\
         I-RDSA [Theorem \ref{th:4}]&Stoch. Non-convex&$m$&$O((d/m)^{1/3}/t^{1/4})$&SZO\\
         \bottomrule
    \end{tabular}
    \caption{Convergence of Frank-Wolfe: Det. refers to deterministic while stoch. refers to stochastic. Memory indicates the number of samples at which the gradients needs to be tracked in the first order case. In the zeroth order case, it indicates the number of directional derivatives being evaluated at one sample. The rates correspond to the rate of decay of $\mathbb{E}[f\left(\mathbf{x}_{t}\right)-f\left(\mathbf{x}^{*}\right)]$ in the convex setting and the Frank-Wolfe duality gap in context of the non-convex setting.}
    \label{tab:FW_convex}
\end{table*}

For the problem in \eqref{eq:prob-def}, it is well known that the primal sub-optimality gap of first order schemes are dimension independent. 
However, algorithms which involve a projection operator might be expensive in practice depending on the structure of $\mathcal{C}$. 
Noting the potentially expensive projection operators, projection free methods such as Frank-Wolfe~\citep{jaggi2013revisiting} have had a resurgence. 
Frank-Wolfe avoids the projection step, and only requires access to a linear minimization oracle, which can be implemented efficiently and needs to be solved to a certain degree of exactness. 
Stochastic versions of Frank-Wolfe have been studied in both the convex~\citep{hazan2012projection,hazan2016variance,mokhtari2018conditional,lan2016conditional} and non-convex~\citep{reddi2016stochastic} setting with access to stochastic first order oracles~(SFO).
However, convergence of stochastic Frank-Wolfe with access to only stochastic zeroth order oracle (SZO) remains unexplored. 

In this paper, we study a setting of the stochastic Frank-Wolfe where a small batch-size~(independent of dimension or the number of iterations) is sampled at each epoch while having access to a zeroth order oracle. 
Unlike, the first order oracle based stochastic Frank-Wolfe, the zeroth order counterpart is only able to generate biased gradient estimates. 
We focus on three different zeroth order gradient approximation schemes, namely, the classical Kiefer Wolfowitz stochastic approximation~(KWSA)~\citep{kiefer1952stochastic}, random directions stochastic approximation~(RDSA)~\citep{nesterov2011random,duchi2015optimal}, and an improvized RDSA~(I-RDSA). 
KWSA samples directional derivatives along the canonical basis directions at each iteration, while RDSA samples one directional derivative at each iteration, and I-RDSA samples $m<d$ directional derivatives at each iteration. 
Na\"ive usage of the biased gradient estimates in the linear minimization step, in addition to the stochasticity of the function evaluations, can lead to potentially diverging iterate sequences. 

To circumvent the potential divergence issue due to non-decaying gradient noise and bias, we use a gradient averaging technique used in \citet{yang2016parallel,ruszczynski2008merit,mokhtari2018conditional} to get a surrogate gradient estimate which reduces the noise and the associated bias. 
The gradient averaging technique intuitively reduces the linear minimization step to that of an inexact minimization if the exact gradient was available. 
For each of the zeroth order optimization schemes, i.e., KWSA, RDSA, and I-RDSA, we derive primal sub-optimality bounds and Frank-Wolfe duality gap bounds and quantify the dependence in terms of the dimension and the number of epochs. 
We show that the primal sub-optimality gap to be of the order $O(d^{1/3}/T^{1/3})$ for RDSA, which improves to $O((d/m)^{1/3}/T^{1/3})$ for I-RDSA, and $O(1/T^{1/3})$ for KWSA at the cost of additional directional derivatives. 
The dimension dependence in zeroth order optimization is unavoidable due to the inherent bias-variance trade-off but nonetheless, the dependence on the number of iterations matches that of its first order counterpart in \citet{mokhtari2018conditional}. Recently in \citep{balasubramanian2018zeroth}, a zeroth order Frank Wolfe algorithm was proposed where the number of gradient directions sampled at each epoch scales linearly with both respect to the number of iterations and dimension of the problem. For the convex case, the number of gradient directions further scales as squared number of iterations. In contrast, we focus on the case where the number of gradient directions sampled at each epoch are independent of the dimension and the number of iterations. Moreover, in \citep{balasubramanian2018zeroth} it is not clear how the primal and dual gap scales with respect to dimension when dimension and iteration independent gradient directions are sampled at each iteration.
Furthermore, we also derive rates for non-convex functions and show the Frank-Wolfe duality gap to be $O(d^{1/3}/T^{1/4})$, where the dependence on the number of iterations matches that of its first order counterpart in \citet{reddi2016stochastic}. To complement the theoretical results, we also demonstrate the efficacy of our algorithm through empirical evaluations on datasets. In particular, we perform experiments on a dataset concerning constrained black box non-convex optimization, where generic first order methods are rendered unusable and show that our proposed algorithm converges to a first order stationary point.

\subsection{Related Work}
\label{subsec:rel_work}
Algorithms for convex optimization with access to a SZO have been studied in \citet{wang2018stochastic,duchi2015optimal,liu2018zeroth,CDCKW2018}, where in \citet{liu2018zeroth} to address constrained optimization a projection step was considered. In the context of projection free methods, \citet{frank1956algorithm} studied the Frank-Wolfe algorithm for smooth convex functions with line search which was extended to encompass inexact linear minimization step in \citet{jaggi2013revisiting}. Subsequently with additional assumptions, the rates for classical Frank-Wolfe was improved in \citet{lacoste2015global,garber2015faster}. Stochastic versions of Frank-Wolfe for convex optimization with number of calls to SFO at each iteration dependent on the number of iterations with additional smoothness assumptions have been studied in \citet{hazan2012projection,hazan2016variance} so as to obtain faster rates, while \citet{mokhtari2018conditional} studied the version with a mini-batch size of $1$. In the context of non-convex optimization, a deterministic Frank-Wolfe algorithm was studied in \citet{lacoste2016convergence}, while \citet{reddi2016stochastic} addressed the stochastic version of Frank-Wolfe and further improved the rates by using variance reduction techniques. Table \ref{tab:FW_convex} gives a summary of the rates of various algorithms.  For the sake of comparison, we do not compare our rates with those of variance reduced versions of stochastic Frank-Wolfe in \citet{reddi2016stochastic,hazan2016variance}, as our proposed algorithm does not employ variance reduction techniques which tend to incorporate multiple restarts and extra memory in order to achieve better rates. However, note that our algorithm can be extended so as to incorporate variance reduction techniques.

\section{Frank-Wolfe: First to Zeroth Order}
\label{sec:fwftoz}
In this paper, the objective is to solve the following optimization problem:
\begin{equation}
\label{eq:obj}
\min_{\mathbf{x}\in\mathcal{C}} f\left(\mathbf{x}\right) = \min_{x\in\mathcal{C}}\mathbb{E}_{\mathbf{y}\sim\mathcal{P}}\left[F\left(\mathbf{x};\mathbf{y}\right)\right],
\end{equation}
where $\mathcal{C}\in\mathbb{R}^{d}$ is a closed convex set, the loss functions and the expected loss functions, $F\left(\cdot;\mathbf{y}\right)$ and $f(\cdot)$ respectively are possibly non-convex. However, in the context of the optimization problem posed in \eqref{eq:obj}, we assume that we have access to a stochastic zeroth order oracle (SZO). On querying a SZO at the iterate $\mathbf{x}_{t}$, yields an unbiased estimate of the loss function $f(\cdot)$ in the form of $F\left(\mathbf{x}_{t};\mathbf{y}_{t}\right)$. Before proceeding to the algorithm and the subsequent results, we revisit preliminaries concerning the Frank-Wolfe algorithm and zeroth order optimization.
\subsection{Background: Frank-Wolfe Algorithm}
The celebrated Frank-Wolfe algorithm is based around approximating the objective by a first-order Taylor approximation. In the case, when exact first order information is available, i.e., one has access to an incremental first order oracle (IFO), a deterministic Frank-Wolfe method involves the following steps: 
\begin{equation}
\begin{aligned}
&\mathbf{v}_{t}=\argmin_{\mathbf{v}\in\mathcal{C}} \langle \nabla f\left(\mathbf{x}_{t}\right), \mathbf{v}\rangle\\
&\mathbf{x}_{t+1}=\left(1-\gamma_{t+1}\right)\mathbf{x}_{t}+\gamma_{t+1}\mathbf{v}_{t},
\end{aligned}
\label{eq:1}
\end{equation}
where $\gamma_{t}=\frac{2}{t+2}$. A linear minimization oracle~(LMO) is queried at every epoch. Note that, the exact minimization in \eqref{eq:1} is a linear program\footnote{Technically speaking, when $\mathcal{C}$ is given by linear constraints.} and can be performed efficiently without much computational overload. It is worth noting that the exact minimization in \eqref{eq:1} can be replaced by an inexact minimization of the following form, where a $\mathbf{v}\in\mathcal{C}$ is chosen to satisfy,
\begin{equation*}
\langle \nabla f\left(\mathbf{x}_{t}\right), \mathbf{v}\rangle\leq\argmin_{\mathbf{v}\in\mathcal{C}} \langle \nabla f\left(\mathbf{x}_{t}\right), \mathbf{v}\rangle+\gamma_{t}C_{1},
\end{equation*}
and the algorithm can be shown to retain the same convergence rate (see, for example \citep{jaggi2013revisiting}).
\subsection{Background: Zeroth Order Optimization}
The crux of zeroth order optimization consists of gradient approximation schemes from appropriately sampled values of the objective function. We briefly describe the few well known zeroth order gradient approximation schemes. The Kiefer-Wolfowitz stochastic approximation~(KWSA, see \citep{kiefer1952stochastic}) scheme approximates the gradient by sampling the objective function along the canonical basis vectors. Formally, gradient estimate can be expressed as:
\begin{equation}
\label{eq:KWSA}
\mathbf{g}(\mathbf{x}_{t};\mathbf{y})=\sum_{i=1}^{d}\frac{F\left(\mathbf{x}_{t}+c_{t}\mathbf{e}_{i};\mathbf{y}\right)-F\left(\mathbf{x}_{t};\mathbf{y}\right)}{c_{t}}\mathbf{e}_{i},
\end{equation}
where $c_t$ is a carefully chosen time-decaying sequence.
KWSA requires $d$ samples at each step to evaluate the gradient. However, in order to avoid sampling the objective function $d$ times, random directions based gradient estimators have been proposed recently~(see, for example \citet{duchi2015optimal,nesterov2011random}). The random directions gradient estimator~(RDSA) involves estimating the directional derivative along a randomly sampled direction from an appropriate probability distribution. Formally, the random directions gradient estimator is given by,
\begin{align}
\label{eq:RD_grad}
\mathbf{g}(\mathbf{x}_{t};\mathbf{y},\mathbf{z}_{t})=\frac{F\left(\mathbf{x}_t+c_{t}\mathbf{z}_{t};\mathbf{y}\right)-F\left(\mathbf{x}_{t};\mathbf{y}\right)}{c_{t}}\mathbf{z}_{t},
\end{align}
where $\mathbf{z}_{t}\in\mathbb{R}^{d}$ is a random vector sampled from a probability distribution such that $\mathbb{E}\left[\mathbf{z}_{t}\mathbf{z}_{t}^{\top}\right]=\mathbf{I}_{d}$ and $c_{t}$ is a carefully chosen time-decaying sequence. With $c_{t}\to 0$, both the gradient estimators in \eqref{eq:KWSA} and \eqref{eq:RD_grad} turn out to be unbiased estimators of the gradient $\nabla f(\mathbf{x}_{t})$.

\section{Zeroth Order Stochastic Frank-Wolfe: Algorithm \& Analysis}
In this section, we start by stating assumptions which are required for our analysis.
\begin{myassump}{A1}
	\label{as:1}
	In problem \eqref{eq:obj}, the set $\mathcal{C}$ is bounded with finite diameter $R$.
\end{myassump}

\begin{myassump}{A2}
	\label{as:2}
	$F$ is convex and Lipschitz continuous with $\sqrt{\E{\left\|\nabla_{x} F(\mathbf{x};\cdot)\right\|^{2}}} \le L_{1}$ for all $\mathbf{x}\in\mathcal{C}$.
\end{myassump}

\begin{myassump}{A3}
	\label{as:3}
	 The expected function $f(\cdot)$ is convex. Moreover, its gradient $\nabla f$ is $L$-Lipschitz continuous over the set $\mathcal{C}$, i.e., for all $x, y \in \mathcal{C}$
	 \begin{equation}
	 \left\|\nabla f(\mathbf{x}) - \nabla f(\mathbf{y})\right\|\leq L\left\|\mathbf{x}-\mathbf{y}\right\|.
	 \end{equation}
	 \end{myassump}

\begin{myassump}{A4}
	\label{as:4}
	The $\mathbf{z}_{t}$'s are drawn from a distribution $\mu$ such that $M(\mu)=\mathbb{E}\left[\left\|\mathbf{z}_{t}\right\|^{6}\right]$ is finite,
	and for any vector $\mathbf{g}\in\mathbb{R}^{d}$, there exists a function $s(d):\mathbb{N}\mapsto\mathbb{R}_{+}$ such that,
	\begin{equation*}
	\mathbb{E}\left[\left\|\langle \mathbf{g},\mathbf{z}_{t}\rangle\mathbf{z}_{t}\right\|^{2}\right] \le s(d)\left\|\mathbf{g}\right\|^{2}.
	\end{equation*}
	\end{myassump}
	
\begin{myassump}{A5}
	\label{as:5}
	The unbiased gradient estimates, $\nabla F\left(\mathbf{x};\mathbf{y}\right)$ of $\nabla f(\mathbf{x})$, i.e., $\mathbb{E}_{\mathbf{y}\sim\mathcal{P}}\left[\nabla F\left(\mathbf{x};\mathbf{y}\right)\right] = \nabla f(\mathbf{x})$ satisfy
	\begin{equation}
	\label{eq:as5}
	\E{\left\|\nabla F\left(\mathbf{x},\mathbf{y}\right)-\nabla f(\mathbf{x})\right\|^{2}} \le \sigma^{2}
	\end{equation}
\end{myassump}

We note that Assumptions \ref{as:1}-\ref{as:3} and \ref{as:5} are standard in the context of stochastic optimization. Assumption \ref{as:4} provides for the requisite moment conditions for the sampling distribution of the directions utilized for finding directional derivatives so as to be able to derive concentration bounds. In particular, if $\mu$ is taken to be uniform on the surface of the $\mathbb{R}^{d}$ Euclidean ball with radius $\sqrt{d}$, then we have that $M(\mu)=d^3$ and $s(d)=d$. Moreover, if $\mu$ is taken to be $\mathcal{N}\left(\mathbf{0},\mathbf{I}_{d}\right)$, then $M(\mu)=d(d+2)(d+4)\approx d^{3}$ and $s(d)=d$. For the rest of the paper, we take $\mu$ to be either uniform on the surface of the $\mathbb{R}^{d}$ Euclidean ball with radius $\sqrt{d}$ or $\mathcal{N}\left(\mathbf{0},\mathbf{I}_{d}\right)$. Before getting into the stochastic case, we demonstrate how a typical zeroth order Frank-Wolfe framework corresponds to an inexact classical Frank-Wolfe optimization in the deterministic setting.
\subsection{Deterministic Zeroth Order Frank-Wolfe} 
The deterministic version of the optimization in \eqref{eq:obj} can be re-stated as follows:
\begin{align}
\label{eq:obj_det}
\min_{\mathbf{x}\in\mathcal{C}} F\left(\mathbf{x}\right).
\end{align}
\begin{algorithm}[tb]
	\caption{Deterministic Zeroth Order Frank Wolfe}
	\label{algo_dfw} 
	\algrenewcommand\algorithmicensure{\textbf{Output:}}
	\begin{algorithmic}[1] 
		\REQUIRE Input, Loss Function $F(x)$, 
		$L$ (Lipschitz constant for the gradients), 
		Convex Set $\mathcal{C}$, 
		Sequences $\gamma_{t}=\frac{2}{t+1}, c_{t}=\frac{L\gamma_{t}}{d}$. 
		\ENSURE: $\mathbf{x}_{T}$ or $\frac{1}{T} \sum_{t=1}^{T}\mathbf{x}_{t}$.
		\STATE Initialize $\mathbf{x}_0\in\mathcal{C}$
		\FOR {$t=0,1,\ldots, T-1$}
			\STATE Compute $\mathbf{g}(\mathbf{x}_{t})=\sum_{i=1}^{d}\frac{F\left(\mathbf{x}_{t}+c_{t}\mathbf{e}_{i}\right)-F\left(\mathbf{x}_{t}\right)}{c_{t}}\mathbf{e}_{i}$,
			\STATE Compute $\mathbf{v}_{t}=\argmin_{\mathbf{s}\in\mathcal{C}} \langle \mathbf{s}, \mathbf{g}(\mathbf{x}_{t})\rangle$,
			\STATE Compute $\mathbf{x}_{t+1} = \left(1-\gamma_{t}\right)\mathbf{x}_{t}+\gamma_{t}\mathbf{v}_{t}$.
			\ENDFOR
	\end{algorithmic}\end{algorithm}
In order to elucidate the equivalence of a typical zeroth order Frank-Wolfe framework corresponds to an inexact classical Frank-Wolfe optimization, we restrict our attention to the Kiefer-Wolfowitz stochastic approximation~(KWSA) for gradient estimation. In particular, the KWSA gradient estimator in \eqref{eq:KWSA} can be expressed as follows:
\begin{align}
\label{eq:grad_KWSA}
&\mathbf{g}(\mathbf{x}_{t})=\sum_{i=1}^{d}\frac{F\left(\mathbf{x}_{t}+c_{t}\mathbf{e}_{i}\right)-F\left(\mathbf{x}_{t}\right)}{c_{t}}\mathbf{e}_{i}\nonumber\\
& = \nabla F(\mathbf{x}_{t})+\sum_{i=1}^{d}\frac{c_{t}}{2}\langle\mathbf{e}_{i}, \nabla^{2} F(\mathbf{x}_{t}+\lambda_{t} c_{t}\mathbf{e}_{i})\mathbf{e}_{i}\rangle \mathbf{e}_{i},
\end{align}
where $\lambda \in [0,1]$. 
The linear optimization step with the current gradient approximation reduces to:
\begin{align}
\label{eq:lmo}
&\langle \mathbf{v}, \mathbf{g}(\mathbf{x}_{t})\rangle = \langle \mathbf{v}, \nabla F(\mathbf{x}_{t})\rangle \nonumber\\&+ \frac{c_{t}}{2}\sum_{i=1}^{d}\langle\mathbf{e}_{i}, \nabla^{2} F(\mathbf{x}_{t}+\lambda_{t} c_{t}\mathbf{e}_{i})\mathbf{e}_{i}\rangle \langle \mathbf{v}, \mathbf{e}_{i}\rangle\nonumber\\
&\Rightarrow\min_{\mathbf{v}\in\mathcal{C}} \langle \mathbf{v}, \mathbf{g}(\mathbf{x}_{t})\rangle \le \min_{\mathbf{s}\in\mathcal{C}}\langle \mathbf{s}, \nabla F(\mathbf{x}_{t})\rangle + \frac{c_{t}LRd}{2}.
\end{align}
In particular, if $c_{t}$ is chosen to be $c_t=\frac{\gamma_{t}}{d}$ and $\gamma_{t}=\frac{2}{t+1}$, we obtain the following bound characterizing the primal gap:
\begin{Theorem}
\label{th:0}
Given the zeroth order Frank-Wolfe algorithm in Algorithm \ref{algo_dfw}, we obtain the following bound:
\begin{equation}
\label{eq:bound1}
F(\bbx_{t}) -F(\bbx^*)  = \frac{Q_{ns}}{t+2},
\end{equation}
where $Q_{ns}=\max\{2(F(\bbx_{0}) -F (\bbx^*)),4LR^{2}\}$.
\end{Theorem}
Theorem \ref{th:0} asserts that with appropriate scaling of $c_t$, i.e., the smoothing parameter for the zeroth order gradient estimator, the iteration dependence of the primal gap matches that of the classical Frank-Wolfe scheme. In particular, for a primal gap of $\epsilon$, the number of iterations needed for the zeroth order scheme in algorithm \ref{algo_dfw} is $O\left(\frac{1}{\epsilon}\right)$, while the number of calls to the linear minimization oracle and zeroth order oracle are given by $O\left(\frac{1}{\epsilon}\right)$ and $O\left(\frac{d}{\epsilon}\right)$ respectively.
The proof of Theorem \ref{th:0} is provided in the appendix \ref{app:1}.

In summary, Theorem \ref{th:0} shows that the deterministic zeroth order Frank-Wolfe algorithm reduces to the inexact classical Frank-Wolfe algorithm with the corresponding primal being dimension independent. However, the dimension independence comes at the cost of querying the zeroth order oracle $d$ times at each iteration. In the sequel, we will focus on the random directions gradient estimator in \eqref{eq:RD_grad} for the stochastic zeroth order Frank-Wolfe algorithm.

\subsection{Zeroth Order Stochastic Frank-Wolfe} 
\label{subsec:zosfw}
In this section, we formally introduce our proposed zeroth order stochastic Frank-Wolfe algorithm. A naive replacement of $\nabla f(\mathbf{x}_{k})$ by its stochastic counterpart, i.e., $\nabla F(\mathbf{x}_{k};\mathbf{y}_{k})$ would make the algorithm potentially divergent due to non-vanishing variance of gradient approximations. Moreover, the naive replacement would lead to the linear minimization constraint to hold only in expectation and thereby potentially also making the algorithm divergent. We use a well known averaging trick to counter this problem which is as follows:
\begin{equation}
\label{eq:2}
\mathbf{d}_{t} = \left(1-\rho_{t}\right)\mathbf{d}_{t-1}+\rho_{t}\mathbf{g}\left(\mathbf{x}_{t},\mathbf{y}_{t}\right),
\end{equation}
where $g\left(\mathbf{x}_{t},\mathbf{y}_{t}\right)$ is a gradient approximation, $\mathbf{d}_0 =\mathbf{0}$ and $\rho_{t}$ is a time-decaying sequence. Technically speaking, such a scheme allows for $\mathbb{E}\left[\left\|\mathbf{d}_{t}-\nabla f\left(\mathbf{x}_{t}\right)\right\|^{2}\right]$ to go to zero asymptotically. With the above averaging scheme, we replace the linear minimization and the subsequent steps as follows:
\begin{align}
\label{eq:3}
&\mathbf{d}_{t} = \left(1-\rho_{t}\right)\mathbf{d}_{t-1}+\rho_{t}\mathbf{g}\left(\mathbf{x}_{t},\mathbf{y}_{t}\right)\nonumber\\
&\mathbf{v}_{t}=\argmin_{\mathbf{v}\in\mathcal{C}} \langle \mathbf{d}_{t}, \mathbf{v}\rangle\nonumber\\
&\mathbf{x}_{t+1}=\left(1-\gamma_{t+1}\right)\mathbf{x}_{t}+\gamma_{t+1}\mathbf{v}_{t}.
\end{align}

We resort to three different gradient approximation schemes for approximating $g\left(\mathbf{x}_{t},\mathbf{y}_{t}\right)$. In particular, in addition to the the KWSA scheme and the random directions scheme, as outlined in \eqref{eq:KWSA} and \eqref{eq:RD_grad}, we employ an improvised random directions gradient estimator~(I-RDSA) by sampling $m$ directions at each time followed by averaging, i.e., $\{\mathbf{z}_{i,t}\}_{i=1}^{m}$ for which we have,
\begin{align}
\label{eq:KW_grad_1}
&\mathbf{g}_{m}(\mathbf{x}_{t};\mathbf{y}_{t},\mathbf{z}_{i,t})\nonumber\\&=\frac{1}{m}\sum_{i=1}^{m}\left(\frac{F\left(\mathbf{x}_t+c_{t}\mathbf{z}_{i,t};\mathbf{y}\right)-F\left(\mathbf{x}_{t};\mathbf{y}\right)}{c_{t}}\mathbf{z}_{i,t}\right).
\end{align}
It is to be noted that the above gradient approximation scheme uses more exactly one data point while utilizing $m$ directional derivatives. In order to quantify the benefits of using such a scheme, we present the statistics concerning the gradient approximation of RDSA and I-RDSA.
We have from \citep{duchi2015optimal} for RDSA,
\begin{align}
\label{eq:grad_approx}
&\mathbb{E}_{\mathbf{z}_{t}\sim\mu, \mathbf{y}_{t}\sim\mathcal{P}}\left[\mathbf{g}(\mathbf{x};\mathbf{y}_{t},\mathbf{z}_{t})\right] = \nabla f\left(\mathbf{x}\right) + c_{t}L\mathbf{v}\left(\mathbf{x},c_{t}\right)\nonumber\\
&\mathbb{E}_{\mathbf{z}_{t}\sim\mu, \mathbf{y}_{t}\sim\mathcal{P}}\left[\left\|\mathbf{g}(\mathbf{x};\mathbf{y}_{t},\mathbf{z}_{t})\right\|^{2}\right] \le 2s(d)\mathbb{E}\left[\left\|\nabla F(\mathbf{x};\mathbf{y}_{t})\right\|^{2}\right]\nonumber\\&+\frac{c_{t}^{2}}{2}L^{2}M(\mu),
\end{align}
Using \eqref{eq:grad_approx}, similar statistics for the improvised RDSA gradient estimator can be evaluated as follows:
\begin{align}
\label{eq:grad_approx_1}
&\mathbb{E}_{\mathbf{z}_{t}\sim\mu, \mathbf{y}_{t}\sim\mathcal{P}}\left[\mathbf{g}_{m}(\mathbf{x};\mathbf{y}_{t},\mathbf{z}_{t})\right] = \nabla f\left(\mathbf{x}\right) + \frac{c_{t}}{m}L\mathbf{v}\left(\mathbf{x},c_{t}\right)\nonumber\\
&\mathbb{E}_{\mathbf{z}_{t}\sim\mu, \mathbf{y}_{t}\sim\mathcal{P}}\left[\left\|\mathbf{g}_{m}(\mathbf{x};\mathbf{y}_{t},\mathbf{z}_{t})\right\|^{2}\right]\le \left(\frac{1+m}{2m}\right)c_{t}^{2}L^{2}M(\mu) \nonumber\\&+2\left(1+\frac{s(d)}{m}\right)\mathbb{E}\left[\left\|\nabla F(\mathbf{x};\mathbf{y}_{t})\right\|^{2}\right],
\end{align}
where $\left\|\mathbf{v}\left(\mathbf{x},c_{t}\right)\right\| \le \frac{1}{2}\mathbb{E}\left[\left\|\mathbf{z}\right\|^{3}\right]$. A proof for \eqref{eq:grad_approx_1} can be found in \citep{liu2018zeroth}. As we will see later the I-RDSA scheme improves the dimension dependence of the primal gap, but it comes at the cost of $m$ calls to the SZO. We are now ready to state the zeroth order stochastic Frank-Wolfe algorithm which is presented in algorithm \ref{algo_sfw}.
\begin{algorithm}[tb]
	\caption{Stochastic Gradient Free Frank Wolfe}
	\label{algo_sfw} 
	\algrenewcommand\algorithmicensure{\textbf{Output:}}
	\begin{algorithmic}[1] 
	\REQUIRE Input, Loss Function $F(x)$, 
	Convex Set $\mathcal{C}$, number of directions $m$, sequences $\gamma_{t}=\frac{2}{t+8}$,
	\begin{itemize}[leftmargin=0mm, itemsep=0mm, partopsep=0pt,parsep=0pt]
	\item[$$] $(\rho_{t},c_{t})_{RDSA} = \left(\frac{4}{d^{1/3}(t+8)^{2/3}},\frac{2}{d^{3/2}(t+8)^{1/3}}\right)$
	\item[$$] $(\rho_{t},c_{t})_{I-RDSA} = \left(\frac{4}{\left(1+\frac{d}{m}\right)^{1/3}(t+8)^{2/3}},\frac{2\sqrt{m}}{d^{3/2}(t+8)^{1/3}}\right)$
	\item[$$] $(\rho_{t},c_{t})_{KWSA} = \left(\frac{4}{(t+8)^{2/3}},\frac{2}{d^{1/2}(t+8)^{1/3}}\right)$. 
	\end{itemize}
	\ENSURE $\mathbf{x}_{T}$ or $\frac{1}{T} \sum_{t=0}^{T-1}\mathbf{x}_{t}$.
	\STATE Initialize $\mathbf{x}_0\in\mathcal{C}$
	\FOR {$t=0,2,\ldots, T-1$}
	\STATE Compute 
	\begin{itemize}[leftmargin=4mm, itemsep=0mm, partopsep=0pt,parsep=0pt]
	\item[$$] KWSA:
	\item[$$] $\mathbf{g}(\mathbf{x}_{t};\mathbf{y})=\sum_{i=1}^{d}\frac{F\left(\mathbf{x}_{t}+c_{t}\mathbf{e}_{i};\mathbf{y}\right)-F\left(\mathbf{x}_{t};\mathbf{y}\right)}{c_{t}}\mathbf{e}_{i}$
	\item[$$] RDSA: Sample $\mathbf{z}_{t}\sim\mathcal{N}(0,\mathbf{I}_{d})$, 
	\item[$$] $\mathbf{g}(\mathbf{x}_{t};\mathbf{y},\mathbf{z}_{t})=\frac{F\left(\mathbf{x}_t+c_{t}\mathbf{z}_{t};\mathbf{y}\right)-F\left(\mathbf{x}_{t};\mathbf{y}\right)}{c_{t}}\mathbf{z}_{t}$
	\item[$$] I-RDSA: Sample $\{\mathbf{z}_{i,t}\}_{i=1}^{m}\sim\mathcal{N}(0,\mathbf{I}_{d})$,
	\item[$$] $\mathbf{g}(\mathbf{x}_{t};\mathbf{y},\mathbf{z}_{t})=\frac{1}{m}\sum_{i=1}^{m}\frac{F\left(\mathbf{x}_t+c_{t}\mathbf{z}_{i,t};\mathbf{y}\right)-F\left(\mathbf{x}_{t};\mathbf{y}\right)}{c_{t}}\mathbf{z}_{i,t}$
	\end{itemize}
	\STATE Compute $\mathbf{d}_{t} = \left(1-\rho_{t}\right)\mathbf{d}_{t-1}+\rho_{t}\mathbf{g}\left(\mathbf{x}_{t},\mathbf{y}_{t}\right)$
	\STATE Compute $\mathbf{v}_{t}=\argmin_{\mathbf{s}\in\mathcal{C}} \langle \mathbf{s}, \mathbf{d}_{t})\rangle$,
	\STATE Compute $\mathbf{x}_{t+1} = \left(1-\gamma_{t}\right)\mathbf{x}_{t}+\gamma_{t}\mathbf{v}_{t}$.
	\ENDFOR
\end{algorithmic}\end{algorithm}
Before the main results, we first study the evolution of the gradient estimates in \eqref{eq:2} and the associated mean square error. 
The following Lemma studies the error of the process $\{\mathbf{d}_{t}\}$ as defined in \eqref{eq:2}.
\begin{Lemma}
\label{le:1}
Let Assumptions \ref{as:1}-\ref{as:5} hold. Given the recursion in \eqref{eq:2}, we have that $\|\nabla f(\bbx_t) - \bbd_t\|^2$ satisfies
\begin{itemize}[leftmargin=5mm, itemsep=0mm, partopsep=0pt,parsep=0pt]
    \item[1)] for the RDSA gradient approximation scheme 
    \begin{align}
\label{eq:le1}
&\E{\|\nabla f(\bbx_t) - \bbd_t\|^2}\leq 2\rho_t^2\sigma^{2}+4\rho_t^2L_{1}^{2}+8\rho_t^2s(d)L_1^{2}\nonumber\\&+2\rho_t^2c_{t}^{2}L^{2}M(\mu)+\frac{2L^2R^2\gamma_{t}^2 }{\rho_t}+\frac{\rho_t}{2}c_{t}^{2}L^{2}M(\mu)\nonumber\\&+\left(1-\frac{\rho_t}{2}\right)\E{\|\nabla f(\bbx_{t-1}) - \bbd_{t-1} \|^2},
\end{align}
\item[2)] for the I-RDSA gradient approximation scheme \begin{align}
\label{eq:le2}
&\E{\|\nabla f(\bbx_t) - \bbd_t\|^2}\leq 2\rho_t^2\left(\sigma^{2}+2L_{1}^{2}\right)\nonumber\\&+\frac{\rho_t}{2m^{2}}c_{t}^{2}L^{2}M(\mu)+8\rho_t^2\left(1+\frac{s(d)}{m}\right)L_1^{2}\nonumber\\
&+\left(\frac{1+m}{2m}\right)\rho_t^2c_{t}^{2}L^{2}M(\mu)+\frac{2L^2R^2\gamma_{t}^2 }{\rho_t}\nonumber\\&+\left(1-\frac{\rho_t}{2}\right)\E{\|\nabla f(\bbx_{t-1}) - \bbd_{t-1} \|^2}
\end{align}
\item[3)] for the KWSA gradient approximation scheme 
\begin{align}
\label{eq:le3}
&\E{\|\nabla f(\bbx_t) - \bbd_t\|^2}\leq 2\rho_t^2\sigma^{2}+2\rho_{t}c_{t}^{2}dL^2\nonumber\\&+\frac{2L^2R^2\gamma_{t}^2 }{\rho_t}+\left(1-\frac{\rho_t}{2}\right)\E{\|\nabla f(\bbx_{t-1}) - \bbd_{t-1} \|^2}.
\end{align}
\end{itemize}
\end{Lemma}
We use the following Lemma so as to study the dynamics of the primal gap.
\begin{Lemma} \label{lemma:bound_on_suboptimality_main}
Consider the zeroth order Frank Wolfe Algorithm in \ref{algo_dfw}. Let Assumptions \ref{as:1}-\ref{as:5} hold. Then, the primal gap $F(\bbx_{t+1})-F(\bbx^*)$ satisfies
	\begin{align}\label{eq:bound_on_suboptimality_main}
	&{F(\bbx_{t+1}) -F(\bbx^*)} 
	\leq (1-\gamma_{t+1}){(F(\bbx_{t}) -F(\bbx^*))}\nonumber\\&+\gamma_{t+1} R\|\nabla F(\bbx_{t})-\bbd_t\| +\frac{LR^2\gamma_{t+1}^2}{2}.
	\end{align}
\end{Lemma}
\begin{proof}
The proof is relegated to the Appendix \ref{app:1}.
\end{proof}
With the above recursions in place, we can now characterize the finite time rates of the mean square errors for the different error approximation schemes. In particular, using Lemma \ref{le:1}, we first state the main result concerning the setting, where the objective is convex. 
\subsubsection{Main Results: Convex Case}
In this section, we state the main results. We first state the main results concerning the primal gap of the proposed algorithm.

\paragraph{Primal Gap:}
We state the main results involving the different gradient approximation schemes for the primal gap, which provide a characterization of $\E{f(\bbx_{t}) -f(\bbx^*) }$.
\begin{Theorem}
	\label{th:1}
	Let Assumptions \ref{as:1}-\ref{as:5} hold. Let the sequence $\gamma_t$ be given by $\gamma_t = \frac{2}{t+8}$.
	\begin{itemize}[leftmargin=5mm, itemsep=0mm, partopsep=0pt,parsep=0pt]
	\item[1)]Then, we have the following primal sub-optimality gap for the algorithm in \ref{algo_sfw}, with the RDSA gradient approximation scheme:
	\begin{align}
	\label{eq:th1}
	\E{f(\bbx_{t}) -f(\bbx^*) } =O\left(\frac{d^{1/3}}{(t+9)^{1/3}}\right).
	\end{align}
	\item[2)]In case of the I-RDSA the gradient approximation scheme, the primal sub-optimality gap is given by,
	\begin{align}
	\label{eq:th2}
	\E{f(\bbx_{t}) -f(\bbx^*) } =O\left(\frac{(d/m)^{1/3}}{(t+9)^{1/3}}\right).
	\end{align}
	\item[3)] Finally, for the KWSA gradient approximation scheme, the primal sub-optimality gap is given by, 
	\begin{align}
    \label{eq:th3}
    \E{f(\bbx_{t}) -f(\bbx^*) } = O\left(\frac{1}{(t+9)^{1/3}}\right).
    \end{align}
\end{itemize}	
\end{Theorem}
Theorem \ref{th:1} quantifies the dimension dependence of the primal gap to be $d^{1/3}$. At the same time the dependence on iterations, i.e., $O(T^{-1/3})$ matches that of the stochastic Frank-Wolfe which has access to first order information as in \citep{mokhtari2018conditional}. The improvement of the rates for I-RDSA and KWSA are at the cost of extra directional derivatives at each iteration. The number of queries to the SZO so as to obtain a primal gap of $\epsilon$, i.e., $\E{f(\bbx_{t}) -f(\bbx^*) } \le \epsilon$ is given by $O\left(\frac{d}{\epsilon^3}\right)$, where the dimension dependence is consistent with zeroth order schemes and cannot be improved on as illustrated in \citep{duchi2015optimal}. The proofs for parts (1), (2) and (3) for theorem \ref{th:1} can be found in the appendix in the sections \ref{app:2}, \ref{app:3} and \ref{app:4} respectively.

\paragraph{Dual Gap:} 
We state the main results involving the different gradient approximation schemes for the dual gap, which provide a characterization of $\mathcal{G}\left(\mathbf{x}\right)=\max_{\mathbf{v}\in\mathcal{C}}\langle \nabla F(\mathbf{x}), \mathbf{x}-\mathbf{v}\rangle$.
\begin{Theorem}
	\label{th:2}
	Let Assumptions \ref{as:1}-\ref{as:5} hold. Let the sequence $\gamma_t$ be given by $\gamma_t = \frac{2}{t+8}$.
	\begin{itemize}[leftmargin=5mm, itemsep=0mm, partopsep=0pt,parsep=0pt]
		\item[1)]Then, we have the following dual gap for the algorithm in \ref{algo_sfw}, with the RDSA gradient approximation scheme:
		\begin{align}
		\label{eq:th11}
		&\E{\min_{t=0,\cdots,T-1}\mathcal{G}\left(\mathbf{x}_{t}\right)} \le \frac{7(F(\mathbf{x}_{0})-F(\mathbf{x}^{\ast}))}{2T}\nonumber\\&+\frac{LR^{2}\ln(T+7)}{T}+\frac{Q^{'}+R\sqrt{2Q}}{2T}(T+7)^{2/3},
		\end{align}
		where $Q=32d^{-1/3}\sigma^2+64d^{-1/3}L_{1}^{2}+128d^{2/3}L_{1}^{2}+2L^2R^2d^{2/3}+416d^{2/3}L^{2}$ and $Q^{'}=\max\{2(f(\bbx_{0}) -f(\bbx^*)),2R\sqrt{2Q}+LR^2/2\}$.
		\item[2)]In case of the I-RDSA the gradient approximation scheme, the dual gap is given by,
		\begin{align}
		\label{eq:th21}
		&\E{\min_{t=0,\cdots,T-1}\mathcal{G}\left(\mathbf{x}_{t}\right)}\le \frac{7(F(\mathbf{x}_{0})-F(\mathbf{x}^{\ast}))}{2T}\nonumber\\&+\frac{LR^{2}\ln(T+7)}{T}+\frac{Q_{ir}^{'}+R\sqrt{2Q_{ir}}}{2T}(T+7)^{2/3},
		\end{align}
		where $Q_{ir}=32\left(1+d/m\right)^{-1/3}\sigma^2+128\left(1+d/m\right)^{2/3}L_{1}^{2}+64\left(1+d/m\right)^{-1/3}L_{1}^{2}+2L^2R^2\left(1+d/m\right)^{2/3}+416\left(1+d/m\right)^{2/3}L^{2}$ and $Q_{ir}^{'}=\max\{2(f(\bbx_{0}) -f(\bbx^*)),2R\sqrt{2Q_{ir}}+LR^2/2$.
		\item[3)] Finally, for the KWSA gradient approximation scheme, the dual gap is given by, 
		\begin{align}
		\label{eq:th31}
		&\E{\min_{t=0,\cdots,T-1}\mathcal{G}\left(\mathbf{x}_{t}\right)} \le \frac{7(F(\mathbf{x}_{0})-F(\mathbf{x}^{\ast}))}{2T}\nonumber\\&+\frac{LR^{2}\ln(T+7)}{T}+\frac{Q_{kw}^{'}+R\sqrt{2Q_{kw}}}{2T}(T+7)^{2/3},
		\end{align}
		where $Q_{kw}=\max \big\{4\|\nabla f(\bbx_0) - \bbd_0\|^2, \; 32\sigma^2+32L^{2}+2L^2R^2\big\}$ and $Q_{kw}^{'}=\max\{2(f(\bbx_{0}) -f(\bbx^*)),2R\sqrt{Q_{kw}}+LR^2/2\}$.
	\end{itemize}	
\end{Theorem}
Theorem \ref{th:2} quantifies the dimension dependence of the Frank-Wolfe duality gap to be $d^{1/3}$. At the same time the dependence on iterations, i.e., $O(T^{-1/3})$ matches that of the primal gap and hence follows that the number of queries to the SZO so as to obtain a Frank-Wolfe duality gap of $\epsilon$, i.e., $\E{\min_{t=0,\cdots,T-1}\mathcal{G}\left(\mathbf{x}_{t}\right)} \le \epsilon$ is given by $O\left(\frac{d}{\epsilon^3}\right)$. In particular, theorem \ref{th:2} asserts that the initial conditions are forgotten as $O(1/T)$. The proofs for parts (1), (2) and (3) for theorem \ref{th:2} can be found in the appendix in the sections \ref{app:2}, \ref{app:3} and \ref{app:4} respectively.

\subsubsection{Zeroth-Order Frank-Wolfe Non-Convex}
We employ the following algorithm for the non-convex stochastic Frank-Wolfe:
\begin{algorithm}[h]
	\caption{Stochastic Gradient Free Frank-Wolfe}
	\label{algo_ncsfw} 
	\algrenewcommand\algorithmicensure{\textbf{Output:}}
	\begin{algorithmic}[1] 
		\REQUIRE Input, Loss Function $F(x)$, Convex Set $\mathcal{C}$, number of directions $m$. Sequences $\gamma=\frac{1}{T^{3/4}}$,
		\begin{itemize}[leftmargin=0mm, itemsep=0mm, partopsep=0pt,parsep=0pt]
				\item[$ $] $(\rho_{t},c_{t}) = \left(\frac{4}{\left(1+\frac{d}{m}\right)^{1/3}(t+8)^{2/3}},\frac{2\sqrt{m}}{(d^{3/2}(t+8)^{1/3})}\right)$
		\end{itemize}
		\ENSURE $\mathbf{x}_{T}$.
		\STATE Initialize $\mathbf{x}_0\in\mathcal{C}$
		\FOR {$t=0,1,\ldots, T-1$}
			\STATE Compute \begin{itemize}
				\item[$$] Sample $\{\mathbf{z}_{i,t}\}_{i=1}^{m}\sim\mathcal{N}(0,\mathbf{I}_{d})$, $\mathbf{g}(\mathbf{x}_{t};\mathbf{y},\mathbf{z}_{t})=\frac{1}{m}\sum_{i=1}^{m}\frac{F\left(\mathbf{x}_t+c_{t}\mathbf{z}_{i,t};\mathbf{y}\right)-F\left(\mathbf{x}_{t};\mathbf{y}\right)}{c_{t}}\mathbf{z}_{i,t}$
				\end{itemize}
			\STATE Compute $\mathbf{d}_{t} = \left(1-\rho_{t}\right)\mathbf{d}_{t-1}+\rho_{t}g\left(\mathbf{x}_{t},\mathbf{y}\right)$
			\STATE Compute $\mathbf{v}_{t}=\argmin_{\mathbf{s}\in\mathcal{C}} \langle \mathbf{s}, \mathbf{d}_{t})\rangle$,
			\STATE Compute $\mathbf{x}_{t+1} = \left(1-\gamma\right)\mathbf{x}_{t}+\gamma\mathbf{v}_{t}$.
		\ENDFOR
	\end{algorithmic}
\end{algorithm}
	
	We use the following assumption concerning the smoothness of the non-convex loss function.
	\begin{myassump}{A6}
		\label{as:6}
		The gradients $\nabla f$ are $L$-Lipschitz continuous over the set $\mathcal{C}$, i.e., for all $x, y \in \mathcal{C}$
		\begin{align*}
		\left\|\nabla f(\mathbf{x}) - \nabla f(\mathbf{y})\right\|\leq L\left\|\mathbf{x}-\mathbf{y}\right\|.
		\end{align*}
	\end{myassump}
	\begin{Theorem}
		\label{th:4}
		Let Assumptions \ref{as:4}-\ref{as:6} hold. Then, we have the following dual gap for iterations $t=0,1,\cdots,T-1$ for the algorithm as described in \eqref{eq:3}
		\begin{align}
		\label{eq:th4}
		\E{\min_{t=0,\cdots,T-1}\mathcal{G}\left(\mathbf{x}_{t}\right)} \leq \frac{Q'}{T^{1/4}}=O\left(\frac{(d/m)^{1/3}}{T^{1/4}}\right),
		\end{align}
		where $Q'=\max\{9^{1/3}(f(\bbx_{0}) -f(\bbx^*)),Q_{nc}R(d/m)^{1/3}\}$. 
	\end{Theorem}
Theorem \ref{th:4} quantifies the dimension dependence of the Frank-Wolfe duality gap for non-convex functions to be $d^{1/3}$. At the same time the dependence on iterations, i.e., $O(T^{-1/4})$ matches that of the rate of SFW in \citep{reddi2016stochastic} and hence follows that the number of queries to the SZO so as to obtain a Frank-Wolfe duality gap of $\epsilon$, i.e., $\E{\min_{t=0,\cdots,T-1}\mathcal{G}\left(\mathbf{x}_{t}\right)} \le \epsilon$ is given by $O\left(\frac{d^{4/3}}{\epsilon^4}\right)$. The proof is relegated to the appendix in section \ref{app:5}.	

\section{Experiments}
\label{sec:expts}
We now present empirical results for zeroth order Frank-Wolfe optimization with an aim to highlight three aspects of our method:
(i) it is accurate even in stochastic case (Section~\ref{sec:exp-stoc})
(ii) it scales to relatively high dimensions (Section~\ref{sec:exp-highd})
(iii) it reaches stationary point in non-convex setting (Section~\ref{sec:exp-nonconvex}).

\paragraph{Methods and Evaluation}
We look at the optimality gap $|f(x_{\text{optimizer}})-f(x^*)|$ as the evaluation metric, where $x_{\text{optimizer}}$ denotes the solution obtained from the employed optimizer and $x^*$ corresponds to true solution.
Most existing zero order optimization techniques like Nelder-Mead simplex \citep{nelder1965simplex} or bound optimization by quadratic approximation (BOBYQA; \citealt{powell2009bobyqa}) can only handle bound constraints, but not arbitrary convex constraints as our method can.
Thus, for all experiments, we could compare proposed zeroth order stochastic Frank-Wolf (0-FW) only with COBYLA, a constrained optimizer by linear approximation, which is popular in engineering fields \citep{powell1994direct}.
For experiments where SFO is available, we additionally compare with stochastic proximal gradient descent (PGD) and first order stochastic Frank-Wolfe method (1-FW).

\begin{figure*}[t]
    \centering
    \vspace{-2mm}
    \begin{subfigure}[t]{0.3\textwidth}
    \includegraphics[width=\linewidth]{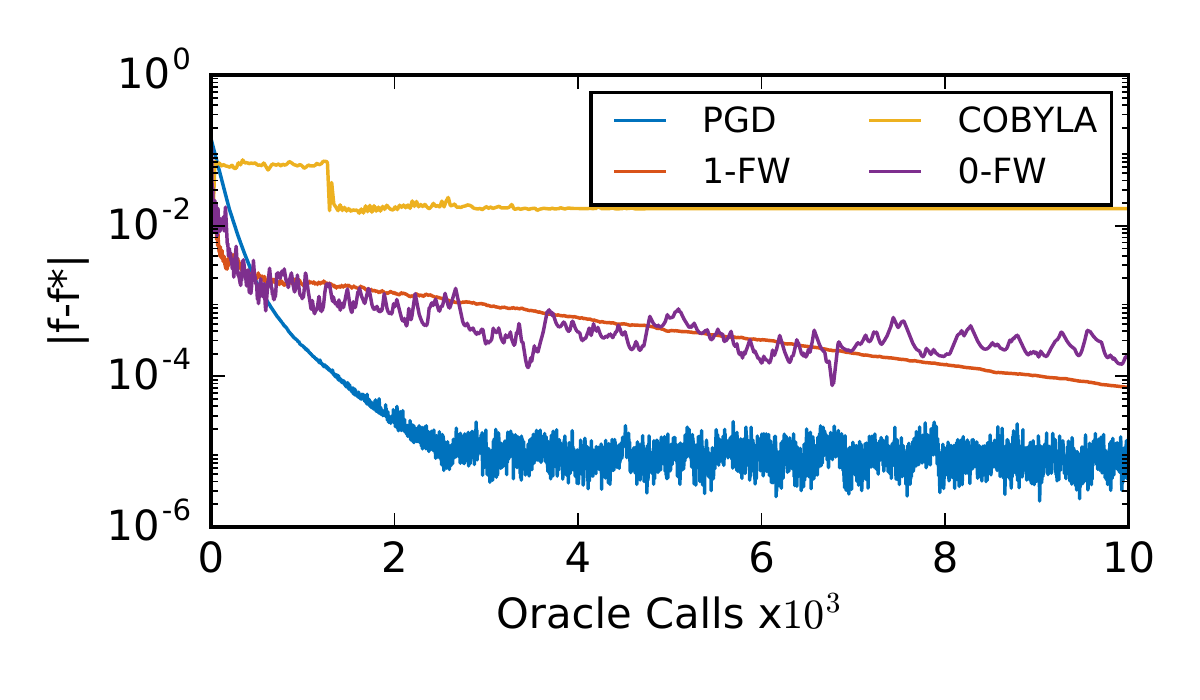}
    \vspace{-4mm}
    \caption{Covtype Dataset}
    \label{fig:covtype}
    \end{subfigure}
    \begin{subfigure}[t]{0.3\textwidth}
    \includegraphics[width=\linewidth]{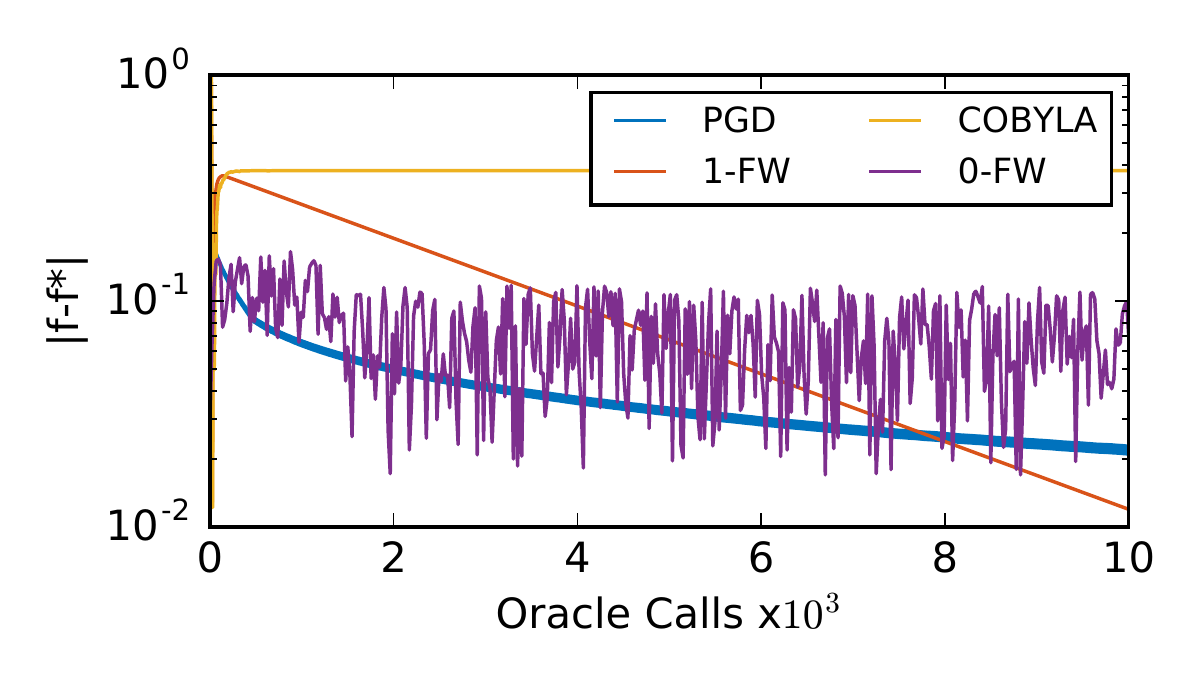}
    \vspace{-4mm}
    \caption{Cox Regression Dataset}
    \label{fig:cox}
    \end{subfigure}
    \begin{subfigure}[t]{0.3\textwidth}
    \includegraphics[width=\linewidth]{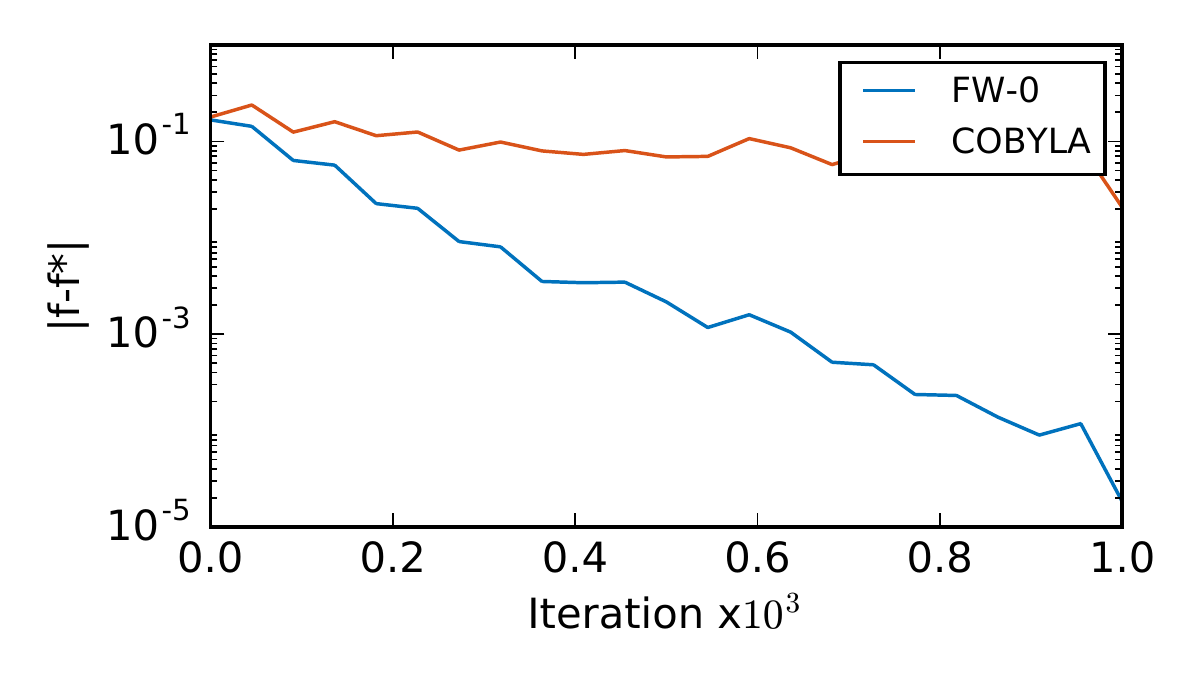}
    \vspace{-4mm}
    \caption{Black Box Optimization Dataset}
    \label{fig:blackbox}
    \end{subfigure}
    \caption{Comparison of proposed zeroth order Frank-Wolfe (0-FW) with first order Frank-Wolfe (1-FW), proximal gradient descent (PGD), and another zero order constrained optimaztion by linear approximation (COBYLA) on various problems.}
    \label{fig:my_label}
\end{figure*}

\subsection{Stochastic Lasso Regression}
\label{sec:exp-stoc}

To study performance of various stochastic optimization, we solve a simple lasso regression on the dataset covtype ($n = 581012$, $d = 54$) from libsvm website\footnote{Available at \url{https://www.csie.ntu.edu.tw/~cjlin/libsvmtools/datasets/}}. 
We use the variant with feature values in $[0,1]$ and solve the following problem:
\begin{equation*}
    \min_{\|\mathbf{w}\|_1 \leq 1} \;\; \frac12 \|\mathbf{y} - \mathbf{X}^\top \mathbf{w}\|_2^2
\end{equation*}
where $\mathbf{X}\in\mathbb{R}^{n \times d}$ represents the feature vectors and $\mathbf{y} \in \mathbb{R}^n$ are the corresponding targets.

For the 0-FW, we used I-RDSA with $m=6$. This problem represents a stochastic setting and from Figure \ref{fig:covtype} we note that the performance of 0-FW matches that of 1-FW in terms of the number of oracle calls to their respective oracles in spite of the dimension involved being $d=54$.

\subsection{High Dimensional Cox Regression}
\label{sec:exp-highd}

To demonstrate efficacy of zeroth order Frank-Wolfe optimization in a moderately high dimensional case, we look at gene expression data. 
In particular, we perform patient survival analysis by solving Cox regression (also known as proportional hazards regression) to relate different gene expression profiles with survival time \citep{sohn2009gradient}. 
We use the Kidney renal clear cell carcinoma dataset\footnote{Available at \url{http://gdac.broadinstitute.org}}, which contains gene expression data for 606 patients (534 with tumor and 72 without tumor) along with survival time information.
We preprocess the dataset by eliminating the rarely expressed genes, i.e. we only keep genes expressed in 50\% of the patients.
This leads to a feature vector $x_i$ of size $9376$ for each patient $i$. 
Also, for each patient $i$, we have the censoring indicator variable $y_i$ that takes the value 0 if patient is alive or 1 if death is observed with $t_i$ denoting the time of death.
In this setup, we can obtain a sparse solution to cox regression by solving the following problem \citep{park2007l1,sohn2009gradient}:
\begin{equation*}
\vspace{-3mm}
    \min_{\|\mathbf{w}\|_1 \leq 10} \frac{1}{n} \sum_{i=1}^n y_i \left\lbrace -\mathbf{x}_i^\top \mathbf{w} + \log\left(\sum_{j\in \mathcal{R}_i} \exp(\mathbf{x}_j^\top \mathbf{w}) \right) \right\rbrace
\end{equation*}
where $\mathcal{R}_i$ is the set of subjects at risk at time $t_i$, i.e. $\mathcal{R}_i = \{j : t_j \geq t_i\}$. 

This problem represents a high-dimensional setting with $d=9376$. For this setup, we take $m=938$ for the I-RDSA scheme of our proposed algorithm. Due to the unavoidable dimension dependence of zeroth order schemes, Figure \ref{fig:cox} shows the gap between 1-FW and 0-FW to be around $2\times$ and thereby reinforcing the result in Theorem \ref{th:1} (2)

\subsection{Black-Box Optimization}
\label{sec:exp-nonconvex}

Finally, we show efficacy of zeroth order Frank-Wolfe optimization in a non-convex setting for a black-box optimization.
Many engineering problems can be posed as optimizing forward models from physics, which are often complicated, do not posses analytical expression, and cannot be differentiated.
We take the example of analyzing electron back-scatter diffraction (EBSD) patterns in order to determine crystal orientation of the sample material.
Such analysis is useful in determining strength, malleability, ductility, etc. of the material along various directions.
Brute-force search has been the primary optimization technique in use \citep{ram2017error}.
For this problem, we use the forward model of EBSD provided by EMSoft\footnote{Software is available at \url{https://github.com/EMsoft-org/EMsoft}}. There are $d=6$ parameters to optimize over the $L_\infty$-ball of radius 1.

This problem represents a non-convex black box optimization setting for which we used $m=1$ for the I-RDSA, i.e. RDSA. Figure \ref{fig:blackbox} shows that our proposed algorithm converges to a first order stationary point there by showing the effectiveness of our proposed algorithm for black-box optimization.

\section{Conclusion}
In this paper, we proposed a stochastic zeroth order Frank-Wolfe algorithm. 
The proposed algorithm does not depend on hard to estimates like Lipschitz constants and thus is easy to deploy in practice.
For the proposed algorithm, we quantified the rates of convergence of the proposed algorithm in terms of the primal gap and the Frank-Wolfe duality gap, which we showed to match its first order counterpart in terms of iterations.
In particular, we showed that the dimension dependence, when one directional derivative is sampled at each iteration to be O($d^{1/3}$). 
We demonstrated the efficacy of our proposed algorithm through experiments on multiple datasets. 
Natural future directions include extending the proposed algorithm to non-smooth functions and incorporating variance reduction techniques to get better rates.


\bibliographystyle{plainnat}
\bibliography{bibliographyJournal,dsprt,glrt,CentralBib}

\newpage
\onecolumn
\appendix

\section{Proofs of Deterministic Frank-Wolfe}
\label{app:1}
\begin{Lemma} \label{lemma:bound_on_suboptimality}
	Consider the proposed zeroth order Frank Wolfe Algorithm. Let Assumptions \ref{as:1}-\ref{as:5} hold. Then, the sub-optimality $F(\bbx_{t+1})-F(\bbx^*)$ satisfies
	\begin{align}\label{eq:bound_on_suboptimality}
	&{F(\bbx_{t+1}) -F(\bbx^*)} 
	\leq (1-\gamma_{t+1}){(F(\bbx_{t}) -F(\bbx^*))}\nonumber\\&+\gamma_{t+1} R\|\nabla F(\bbx_{t})-\bbd_t\| +\frac{LR^2\gamma_{t+1}^2}{2}.
	\end{align}
\end{Lemma}
\begin{proof}
The $L$-smoothness of the function $f$ yields the following upper bound on $f(\bbx_{t+1})$:
{\small\begin{align}\label{eq:proof_lemma_convex_100}
&f(\bbx_{t+1})\leq f(\bbx_{t}) +\nabla f(\bbx_{t})^T (\bbx_{t+1}-\bbx_t) +\frac{L}{2} \|\bbx_{t+1}-\bbx_t\|^2\nonumber\\
&=f(\bbx_{t}) +\gamma_{t+1} (\nabla f(\bbx_{t})-\bbd_t)^T (\bbv_t-\bbx_t) +\gamma_{t+1}\bbd_t^T (\bbv_t-\bbx_t) \nonumber\\&+\frac{L\gamma_{t+1}^2}{2} \|\bbv_{t}-\bbx_t\|^2
\end{align}}
Since $\langle \bbx^*, \bbd_t \rangle \geq \min_{v\in \ccalC} \{ \langle \bbv, \bbd_t \rangle\}  = \langle \bbv_t, \bbd_t \rangle$, we have,
{\small\begin{align}\label{eq:proof_lemma_convex_300}
&f(\bbx_{t+1})\leq f(\bbx_{t})+\gamma_{t+1} (\nabla f(\bbx_{t})-\bbd_t)^T (\bbv_t-\bbx_t)\nonumber\\&+\gamma_{t+1}\bbd_t^T (\bbx^*-\bbx_t) +\frac{L\gamma_{t+1}^2}{2} \|\bbv_{t}-\bbx_t\|^2\nonumber\\
&\leq f(\bbx_{t}) +\gamma_{t+1} (\nabla f(\bbx_{t})-\bbd_t)^T (\bbv_t-\bbx^*) \nonumber\\&+\gamma_{t+1}\nabla f(\bbx_{t})^T (\bbx^*-\bbx_t)+\frac{LR\gamma_{t+1}^2}{2} \|\bbv_{t}-\bbx_t\|^2.
\end{align}}
Using Cauchy-Schwarz inequality, we have,
{\small\begin{align}\label{eq:proof_lemma_convex_500}
&f(\bbx_{t+1})  \leq f(\bbx_{t}) +\gamma_{t+1} \|\nabla f(\bbx_{t})-\bbd_t\| \|\bbv_t-\bbx^*\|\nonumber\\&-\gamma_{t+1}( f(\bbx_{t})-f(\bbx^*))+\frac{L\gamma_{t+1}^2}{2}\|\bbv_t-\bbx^*\| ^2\nonumber\\
& \leq f(\bbx_{t}) +\gamma_{t+1} R\|\nabla f(\bbx_{t})-\bbd_t\| -\gamma_{t+1}( f(\bbx_{t})-f(\bbx^*))\nonumber\\&+\frac{LR^2\gamma_{t+1}^2}{2},
\end{align}}
and subtracting $f(\bbx^*)$ from both sides of \eqref{eq:proof_lemma_convex_500}, we have,
\begin{align}\label{eq:proof_lemma_convex_600}
&f(\bbx_{t+1}) -f(\bbx^*) \leq (1-\gamma_{t+1})(f(\bbx_{t}) -f(\bbx^*))\nonumber\\&+\gamma_{t+1} R\|\nabla f(\bbx_{t})-\bbd_t\| +\frac{LR^2\gamma_{t+1}^2}{2}.
\end{align}
\end{proof}
\begin{proof}[Proof of Theorem \ref{th:0}]
We have, from Lemma \ref{lemma:bound_on_suboptimality},
\begin{align}
\label{eq:bound}
&{F(\bbx_{t+1}) -F(\bbx^*) } \leq (1-\gamma_{t+1}){(F(\bbx_{t}) -F(\bbx^*))}\nonumber\\&+\gamma_{t+1} R\|\nabla F(\bbx_{t})-\mathbf{g}(\mathbf{x}_{t})\| +\frac{LR^2\gamma_{t+1}^2}{2} \nonumber\\
&\Rightarrow {F(\bbx_{t+1}) -F(\bbx^*) } \leq (1-\gamma_{t+1}){(F(\bbx_{t}) -F(\bbx^*))}\nonumber\\&+\frac{c_{t+1}d}{2}\gamma_{t+1} R^{2} +\frac{LR^2\gamma_{t+1}^2}{2}.
\end{align}
From, \eqref{eq:bound}, we have,
\begin{align}
\label{eq:th0proof_1}
{F(\bbx_{t+1}) -F(\bbx^*) } \leq (1-\gamma_{t+1}){(F(\bbx_{t}) -F(\bbx^*))}&+LR^2\gamma_{t+1}^2.
\end{align}
We use Lemma \ref{lemma-estimation} to derive the primal gap which then yields,
\begin{align}
\label{eq:bound1_proof}
F(\bbx_{t}) -F(\bbx^*)  = \frac{Q_{ns}}{t+2},
\end{align}
where $Q_{ns}=\max\{2(F(\bbx_{0}) -F (\bbx^*)),4LR^{2}\}$.
\end{proof}
\section{Proofs of Zeroth Order Stochastic Frank Wolfe: RDSA}
\label{app:2}
\begin{proof}[Proof of Lemma \ref{le:1} (1)]
Use the definition $\bbd_t := (1-\rho_t) \bbd_{t-1} + \rho_t g(\bbx_t; \bby_t,\bbz_t)$  to write the difference $\|\nabla f(\bbx_t) - \bbd_t\|^2$ as
\begin{align}\label{proof:bound_on_grad_100}
&\|\nabla f(\bbx_t) - \bbd_t\|^2  =\|\nabla f(\bbx_t) -(1-\rho_t) \bbd_{t-1}\nonumber\\& - \rho_t g(\bbx_t; \bby_t,\bbz_t)\|^2.
\end{align}
Add and subtract the term $(1-\rho_t)\nabla f(\bbx_{t-1})$ to the right hand side of \eqref{proof:bound_on_grad_100}, regroup the terms and expand the squared term to obtain
{\small\begin{align}\label{proof:bound_on_grad_200}
&\|\nabla f(\bbx_t) - \bbd_t\|^2 \nonumber\\
& =\|\nabla f(\bbx_t)-(1-\rho_t)\nabla f(\bbx_{t-1})+(1-\rho_t)\nabla f(\bbx_{t-1}) \nonumber\\&-(1-\rho_t) \bbd_{t-1} - \rho_t g(\bbx_t; \bby_t,\bbz_t)\|^2 \nonumber\\
&=\rho_t^2\|\nabla f(\bbx_t)-g(\bbx_t; \bby_t,\bbz_t)\|^2  
+(1-\rho_t)^2\|\nabla f(\bbx_{t-1}) - \bbd_{t-1} \|^2\nonumber\\
&+(1-\rho_t)^2\|\nabla f(\bbx_t)-\nabla f(\bbx_{t-1})\|^2\nonumber\\ \nonumber\\
&  \qquad + 2\rho_t (1-\rho_t)(\nabla f(\bbx_t)-g(\bbx_t; \bby_t,\bbz_t))^T(\nabla f(\bbx_t)-\nabla f(\bbx_{t-1})) \nonumber\\
&  \qquad + 2\rho_t (1-\rho_t)(\nabla f(\bbx_t)-g(\bbx_t; \bby_t,\bbz_t))^T(\nabla f(\bbx_{t-1}) - \bbd_{t-1} )\nonumber\\
&  \qquad + 2 (1-\rho_t)^2(\nabla f(\bbx_t)-\nabla f(\bbx_{t-1}))^T(\nabla f(\bbx_{t-1}) - \bbd_{t-1} ).
\end{align}}
Compute the expectation $\E{ (.) \mid\ccalF_t}$ for both sides of \eqref{proof:bound_on_grad_200}, where $\mathcal{F}_{t}$ is the $\sigma$-algebra given by $\{\{\mathbf{y}_{s}\}_{s=0}^{t-1},\{\mathbf{z}_{s}\}_{s=0}^{t-1}\}$to obtain
{\small\begin{align}\label{proof:bound_on_grad_300}
&\E{\|\nabla f(\bbx_t) - \bbd_t\|^2\mid\ccalF_t} \nonumber\\
&=\rho_t^2\E{\|\nabla f(\bbx_t)-g(\bbx_t; \bby_t,\bbz_t)\|^2\mid\ccalF_t}  
\nonumber\\&+(1-\rho_t)^2\|\nabla f(\bbx_t)-\nabla f(\bbx_{t-1})\|^2\nonumber\\
&
+(1-\rho_t)^2\|\nabla f(\bbx_{t-1}) - \bbd_{t-1} \|^2 \nonumber\\&+ 2 (1-\rho_t)^2(\nabla f(\bbx_t)-\nabla f(\bbx_{t-1}))^T(\nabla f(\bbx_{t-1}) - \bbd_{t-1} )\nonumber\\
&+ 2\rho_t (1-\rho_t)\mathbb{E}\left[(\nabla f(\bbx_t)-g(\bbx_t; \bby_t,\bbz_t))^T(\nabla f(\bbx_t)-\nabla f(\bbx_{t-1}))\mid\ccalF_t\right]\nonumber\\&+ 2\rho_t (1-\rho_t)\mathbb{E}\left[(\nabla f(\bbx_t)-g(\bbx_t; \bby_t,\bbz_t))^T(\nabla f(\bbx_{t-1}) - \bbd_{t-1} )\mid\ccalF_t\right]\nonumber\\
&\le\rho_t^2\E{\|\nabla f(\bbx_t)-g(\bbx_t; \bby_t,\bbz_t)\|^2\mid\ccalF_t}  
\nonumber\\&+(1-\rho_t)^2\|\nabla f(\bbx_t)-\nabla f(\bbx_{t-1})\|^2\nonumber\\
&+(1-\rho_t)^2\|\nabla f(\bbx_{t-1}) - \bbd_{t-1} \|^2\nonumber\\&+(1-\rho_t)^2\beta_t\|\nabla f(\bbx_{t-1}) - \bbd_{t-1}\|^2+\frac{(1-\rho_t)^2}{\beta_t}\|\nabla f(\bbx_t)-\nabla f(\bbx_{t-1}) \|^2\nonumber\\
&+2\rho_t (1-\rho_t)(c_{t}L\mathbf{v}\left(\mathbf{x},c_{t}\right))^{\top}(\nabla f(\bbx_t)-\nabla f(\bbx_{t-1}))\nonumber\\
&+2\rho_t (1-\rho_t)(c_{t}L\mathbf{v}\left(\mathbf{x},c_{t}\right))^{\top}(\nabla f(\bbx_{t-1}) - \bbd_{t-1})\nonumber\\
&\le \rho_t^2\E{\|\nabla f(\bbx_t)-g(\bbx_t; \bby_t,\bbz_t)\|^2\mid\ccalF_t}  
\nonumber\\&+(1-\rho_t)^2\|\nabla f(\bbx_t)-\nabla f(\bbx_{t-1})\|^2\nonumber\\
&+(1-\rho_t)^2\|\nabla f(\bbx_{t-1}) - \bbd_{t-1} \|^2 \nonumber\\&+(1-\rho_t)^2\beta_t\|\nabla f(\bbx_{t-1}) - \bbd_{t-1}\|^2+\frac{(1-\rho_t)^2}{\beta_t}\|\nabla f(\bbx_t)-\nabla f(\bbx_{t-1}) \|^2\nonumber\\
&+2\rho_t (1-\rho_t)c_{t}^{2}\left\|L\mathbf{v}\left(\mathbf{x},c_{t}\right)\right\|^{2}+\rho_t (1-\rho_t)\left\|\nabla f(\bbx_t)-\nabla f(\bbx_{t-1})\right\|^{2}\nonumber\\
&+\rho_t (1-\rho_t)\left\|\nabla f(\bbx_{t-1}) - \bbd_{t-1}\right\|^{2}\nonumber\\
&\Rightarrow \E{\|\nabla f(\bbx_t) - \bbd_t\|^2} \nonumber\\&\le \rho_t^2\E{\|\nabla f(\bbx_t)-\nabla F(\bbx_t,\bby_t)+\nabla F(\bbx_t,\bby_t)-g(\bbx_t; \bby_t,\bbz_t)\|^2} \nonumber\\
&+(1-\rho_t)^2\E{\|\nabla f(\bbx_t)-\nabla f(\bbx_{t-1})\|^2}\nonumber\\
&+(1-\rho_t)^2\|\E{\nabla f(\bbx_{t-1}) - \bbd_{t-1} \|^2} \nonumber\\&+(1-\rho_t)^2\beta_t\E{\|\nabla f(\bbx_{t-1}) - \bbd_{t-1}\|^2}\nonumber\\&+\frac{(1-\rho_t)^2}{\beta_t}\E{\|\nabla f(\bbx_t)-\nabla f(\bbx_{t-1}) \|^2}\nonumber\\
&+\frac{\rho_t}{4} (1-\rho_t)c_{t}^{2}L^{2}M(\mu)+\rho_t (1-\rho_t)\E{\left\|\nabla f(\bbx_t)-\nabla f(\bbx_{t-1})\right\|^{2}}\nonumber\\
&+\rho_t (1-\rho_t)\E{\left\|\nabla f(\bbx_{t-1}) - \bbd_{t-1}\right\|^{2}}\nonumber\\
&\le 2\rho_t^2\E{\|\nabla f(\bbx_t)-\nabla F(\bbx_t,\bby_t)\|^{2}}\nonumber\\&+2\rho_{t}^{2}\E{\|\nabla F(\bbx_t,\bby_t)-g(\bbx_t; \bby_t,\bbz_t)\|^2} \nonumber\\
&+\left(1-\rho_t+\frac{(1-\rho_t)^2}{\beta_{t}}\right)\E{\|\nabla f(\bbx_t)-\nabla f(\bbx_{t-1})\|^2}
\nonumber\\&+\left(1-\rho_t+(1-\rho_t)^2\beta_{t}\right)\E{\|\nabla f(\bbx_{t-1}) - \bbd_{t-1} \|^2} \nonumber\\&
+\frac{\rho_t}{2} (1-\rho_t)c_{t}^{2}L^{2}M(\mu)\nonumber\\
&\le 2\rho_t^2\sigma^{2}+4\rho_{t}^{2}\E{\|\nabla F(\bbx_t,\bby_t)\|^{2}}+4\rho_{t}^{2}\E{\|g(\bbx_t; \bby_t,\bbz_t)\|^2} \nonumber\\
&+\left(1-\rho_t+\frac{(1-\rho_t)^2}{\beta_{t}}\right)\E{\|\nabla f(\bbx_t)-\nabla f(\bbx_{t-1})\|^2}
\nonumber\\&+\left(1-\rho_t+(1-\rho_t)^2\beta_{t}\right)\E{\|\nabla f(\bbx_{t-1}) - \bbd_{t-1} \|^2} \nonumber\\&
+\frac{\rho_t}{2} (1-\rho_t)c_{t}^{2}L^{2}M(\mu)\nonumber\\
&\le 2\rho_t^2\sigma^{2}+4\rho_{t}^{2}L_1^2+8\rho_{t}^{2}s(d)L_{1}^{2}+2\rho_{t}^{2}c_{t}^{2}L^{2}M(\mu) \nonumber\\
&+\left(1-\rho_t+\frac{(1-\rho_t)^2}{\beta_{t}}\right)\E{\|\nabla f(\bbx_t)-\nabla f(\bbx_{t-1})\|^2}
\nonumber\\&+\left(1-\rho_t+(1-\rho_t)^2\beta_{t}\right)\E{\|\nabla f(\bbx_{t-1}) - \bbd_{t-1} \|^2} \nonumber\\&
+\frac{\rho_t}{2}c_{t}^{2}L^{2}M(\mu),
\end{align}}
where we used the gradient approximation bounds as stated in \eqref{eq:grad_approx} and used Young's inequality to substitute the inner products and in particular substituted $2\langle \nabla f(\bbx_t)-\nabla f(\bbx_{t-1}) , \nabla f(\bbx_{t-1}) - \bbd_{t-1} \rangle$ by the upper bound $\beta_t \|\nabla f(\bbx_{t-1}) - \bbd_{t-1}\|^2+(1/\beta_t)\|\nabla f(\bbx_t)-\nabla f(\bbx_{t-1}) \|^2$ where $\beta_t>0$ is a free parameter. \\
By assumption \ref{as:4}, the norm $\|\nabla f(\bbx_t)-\nabla f(\bbx_{t-1}) \|$ is bounded above by $L\|\bbx_t-\bbx_{t-1}\|$. In addition, the condition in Assumption \ref{as:1} implies that $L\|\bbx_t-\bbx_{t-1}\|= L\gamma_{t}\|\bbv_t-\bbx_t\| \leq \gamma_{t}LR$. Therefore, we can replace $\|\nabla f(\bbx_t)-\nabla f(\bbx_{t-1}) \|$ by its upper bound $\gamma_{t}LR$ and since we assume that $\rho_t\leq1$ we can replace all the terms $(1-\rho_t)^2$. Furthermore, using $\beta_{t}:=\rho_t/2$ we have,
\begin{align}\label{proof:bound_on_grad_500}
&\E{\|\nabla f(\bbx_t) - \bbd_t\|^2} \nonumber\\
&\leq 2\rho_t^2\sigma^{2}+4\rho_t^2L_{1}^{2}+8\rho_t^2s(d)L_1^{2}+2\rho_t^2c_{t}^{2}L^{2}M(\mu)\nonumber\\
&+\gamma_{t}^2(1-\rho_t)\left(1+\frac{2}{\rho_t}\right)L^2R^2+\frac{\rho_t}{2}c_{t}^{2}L^{2}M(\mu)\nonumber\\
&+(1-\rho_t)\left(1+\frac{\rho_t}{2}\right)\E{\|\nabla f(\bbx_{t-1}) - \bbd_{t-1} \|^2} .
\end{align}
Now using the inequalities $(1-\rho_t) (1+(2/\rho_t))\leq (2/\rho_t)$ and $(1-\rho_t)(1+({\rho_t}/{2}))\leq (1-\rho/2)$ we obtain
\begin{align}\label{proof:bound_on_grad_700}
&\E{\|\nabla f(\bbx_t) - \bbd_t\|^2}\leq 2\rho_t^2\sigma^{2}+4\rho_t^2L_{1}^{2}\nonumber\\&+8\rho_t^2s(d)L_1^{2}+2\rho_t^2c_{t}^{2}L^{2}M(\mu)\nonumber\\
&+\frac{2L^2R^2\gamma_{t}^2 }{\rho_t}+\frac{\rho_t}{2}c_{t}^{2}L^{2}M(\mu)\nonumber\\&+\left(1-\frac{\rho_t}{2}\right))\E{\|\nabla f(\bbx_{t-1}) - \bbd_{t-1} \|^2}.
\end{align}
\end{proof}
Then, we have, from Lemma \ref{lemma:bound_on_suboptimality}
\begin{align}
&\E{f(\bbx_{t+1}) -f(\bbx^*) } 
\leq (1-\gamma_{t+1})\E{(f(\bbx_{t})-f(\bbx^*))}\nonumber\\&+\gamma_{t+1} R\E{\|\nabla f(\bbx_{t})-\bbd_t\|} +\frac{LR^2\gamma_{t+1}^2}{2},
\end{align}
and then by using Jensen's inequality, we obtain,
{\small\begin{align}\label{eq:proof_of_rate_thm_cvx_100}
&\E{f(\bbx_{t+1}) -f(\bbx^*)} 
\leq (1-\gamma_{t+1})\E{(f(\bbx_{t}) -f(\bbx^*))}\nonumber\\&+\gamma_{t+1} R\sqrt{\E{\|\nabla f(\bbx_{t})-\bbd_t\|^2}} +\frac{LR^2\gamma_{t+1}^2}{2}.
\end{align}}

We state a Lemma next which will be crucial for the rest of the paper.
\begin{Lemma}
	\label{lemma-estimation}
	Let $z(k)$ be a non-negative (deterministic) sequence satisfying:
	\[
	z(k+1) \leq (1-r_1(k))\,z_1(k) + r_2(k),
	\]
	where $\{r_1(k)\}$ and $\{r_2(k)\}$
	are deterministic sequences with
	\begin{eqnarray*}
		\frac{a_1}{(k+1)^{\delta_1}}
		\leq r_1(k) \leq 1 \,\,\,\mathrm{and}\,\,\,
		r_2(k) \leq
		\frac{a_2}{(k+1)^{2\delta_1}},
	\end{eqnarray*}
	with $a_1 >0$ , $a_2 >0$ , $1 > \delta_1 > 1/2$ and $k_0 \geq 1$. Then, 
	\begin{align*}
	z(k+1) \leq  \exp\left(-\frac{a_1\delta_1(k+k_0)^{1-\delta_1}}{4(1-\delta_1)}\right)\left(z(0)+\frac{a_{2}}{k_0^{\delta_1}(2\delta_1-1)}\right)+\frac{a_{2}2^{\delta_1}}{a_{1}\left(k+k_0\right)^{\delta_1}}.
	\end{align*}
\end{Lemma}
\begin{proof}[Proof of Lemma \ref{lemma-estimation}]
We have,
\begin{align}
\label{eq:pita-1}
&z(k+1) \le \prod_{l=0}^{k}\left(1-\frac{a_1}{(l+k_0)^{\delta_1}}\right)z(0)\sum_{l=0}^{\lfloor \frac{k}{2}\rfloor -1}\prod_{m=l+1}^{k}\left(1-\frac{a_1}{(m+k_0)^{\delta_1}}\right)\frac{a_2}{(k+k_0)^{2\delta_1}}\nonumber\\&+\sum_{l=\lfloor \frac{k}{2}\rfloor}^{k}\prod_{m=l+1}^{k}\left(1-\frac{a_1}{(m+k_0)^{\delta_1}}\right)\frac{a_2}{(k+k_0)^{2\delta_1}}\nonumber\\
&\le \exp\left(\sum_{l=0}^{k}\left(1-\frac{a_1}{(l+k_0)^{\delta_1}}\right)\right)z(0)+\prod_{m=l+1}^{k}\left(1-\frac{a_1}{(m+k_0)^{\delta_1}}\right)\sum_{l=0}^{\lfloor \frac{k}{2}\rfloor -1}\frac{a_2}{(k+k_0)^{2\delta_1}}\nonumber\\
&+\frac{a_{2}2^{\delta_1}}{a_{1}\left(k+k_0\right)^{\delta_1}}\sum_{l=\lfloor \frac{k}{2}\rfloor}^{k}\prod_{m=l+1}^{k}\left(1-\frac{a_1}{(m+k_0)^{\delta_1}}\right)\frac{a_1}{(k+k_0)^{\delta_1}}\nonumber\\
&\le \exp\left(-\sum_{l=0}^{k}\frac{a_1}{(l+k_0)^{\delta_1}}\right)z(0)+\frac{a_{2}}{a_{1}k_0^{\delta_1}}\exp\left(-\sum_{m=\lfloor \frac{k}{2}\rfloor}^{k}\frac{a_1}{(m+k_0)^{\delta_1}}\right)\sum_{l=0}^{\lfloor \frac{k}{2}\rfloor -1}\frac{a_1}{(k+k_0)^{2\delta_1}}\nonumber\\
&+\frac{a_{2}2^{\delta_1}}{a_{1}\left(k+k_0\right)^{\delta_1}}\sum_{l=\lfloor \frac{k}{2}\rfloor}^{k}\left(\prod_{m=l+1}^{k}\left(1-\frac{a_1}{(m+k_0)^{\delta_1}}\right)-\prod_{m=l}^{k}\left(1-\frac{a_1}{(m+k_0)^{\delta_1}}\right)\right)\nonumber\\
&\le \exp\left(-\sum_{l=0}^{k}\frac{a_1}{(l+k_0)^{\delta_1}}\right)z(0)+\frac{a_{2}2^{\delta_1}}{a_{1}\left(k+k_0\right)^{\delta_1}}+\frac{a_{2}}{a_{1}k_0^{\delta_1}}\exp\left(-\sum_{m=\lfloor \frac{k}{2}\rfloor}^{k}\frac{a_1}{(m+k_0)^{\delta_1}}\right)\sum_{l=0}^{\lfloor \frac{k}{2}\rfloor -1}\frac{a_1}{(k+k_0)^{2\delta_1}}\nonumber\\
&\le \exp\left(-\sum_{l=0}^{k}\frac{a_1}{(l+k_0)^{\delta_1}}\right)z(0)+\frac{a_{2}2^{\delta_1}}{a_{1}\left(k+k_0\right)^{\delta_1}}+\frac{a_{2}}{k_0^{\delta_1}}\exp\left(-\frac{a_1\delta_1}{4(1-\delta_1)}(k+k_0)^{1-\delta_1}\right)\frac{1}{2\delta_1-1},
\end{align}
where we used the inequality that, 
\begin{align*}
&\sum_{m=\lfloor\frac{k}{2}\rfloor}^{k}\frac{1}{(m+k_0)^{\delta_1}} \geq \frac{1}{2(1-\delta_1)}(k+k_0)^{1-\delta_1}-\frac{1}{2(1-\delta_1)}\left(\frac{k}{2}+k_0\right)^{1-\delta_1}\nonumber\\
&\geq \frac{1}{2^{1+\delta_1}(1-\delta_1)}(k+k_0)^{1-\delta_1}\left(2^{1-\delta_1}-1-\frac{(1-\delta_1)k_0}{k+k_0}\right)\geq \frac{\delta_1}{4(1-\delta_1)}(k+k_0)^{1-\delta_1}
\end{align*}
Following up with \eqref{eq:pita-1}, we have,
\begin{align}
\label{eq:pita-2}
&z(k+1) \le \exp\left(-\sum_{l=0}^{k}-\frac{a_1}{(l+k_0)^{\delta_1}}\right)z(0)+\frac{a_{2}2^{\delta_1}}{a_{1}\left(k+k_0\right)^{\delta_1}}+\frac{a_{2}}{k_0^{\delta_1}}\exp\left(-\frac{a_1\delta_1}{4(1-\delta_1)}(k+k_0)^{1-\delta_1}\right)\frac{1}{2\delta_1-1}\nonumber\\
&\le \exp\left(-\frac{a_1\delta_1(k+k_0)^{1-\delta_1}}{4(1-\delta_1)}\right)\left(z(0)+\frac{a_{2}}{k_0^{\delta_1}(2\delta_1-1)}\right)+\frac{a_{2}2^{\delta_1}}{a_{1}\left(k+k_0\right)^{\delta_1}}.
\end{align}
For $\delta=2/3$, we have,
\begin{align*}
z(k+1) \le \exp\left(-\frac{a_1(k+k_0)^{1/3}}{2}\right)\left(z(0)+\frac{3a_{2}}{k_0^{2/3}}\right)+\frac{a_{2}2^{2/3}}{a_{1}\left(k+k_0\right)^{2/3}}.
\end{align*}
\end{proof}
\begin{proof}[Proof of Theorem \ref{th:1} (1)]
Now using the result in Lemma \ref{lemma-estimation} we can characterize the convergence of the sequence of expected errors $\E{\|\nabla f(\bbx_{t})-\bbd_t\|^2}$ to zero. To be more precise, using the result in Lemma \ref{le:1} and setting $\gamma_t=2/(t+8)$, $\rho_t=4/d^{1/3}(t+8)^{2/3}$ and $c_{t}=2/\sqrt{M(\mu)}(t+8)^{1/3}$ for any $\epsilon > 0$ to obtain
{\small\begin{align}\label{eq:almost_done}
& \E{\|\nabla f(\bbx_t) - \bbd_t\|^2}\nonumber\\
&\leq \left(1-\frac{2}{d^{1/3}(t+8)^{2/3}}\right)\E{\|\nabla F(\bbx_{t-1}) - \bbd_{t-1}\|^2}\nonumber\\&+\frac{32d^{-1/3}\sigma^2+64d^{-1/3}L_{1}^{2}+128d^{2/3}L_{1}^{2}+2L^2R^2d^{2/3}+416d^{2/3}L^{2} }{(t+8)^{4/3}}.
\end{align}}
 According to the result in Lemma \ref{lemma-estimation}, the inequality in  \eqref{eq:almost_done} implies that 
 \begin{equation}\label{eq:done}
 \E{\|\nabla f(\bbx_t) - \bbd_t\|^2} \leq \overline{Q}+\frac{Q}{(t+8)^{2/3}} \leq \frac{2Q}{(t+8)^{2/3}},
 \end{equation}
 where $Q=32d^{-1/3}\sigma^2+64d^{-1/3}L_{1}^{2}+128d^{2/3}L_{1}^{2}+2L^2R^2d^{2/3}+416d^{2/3}L^{2}$, where $\overline{Q}$ is a function of $\E{\|\nabla f(\bbx_0) - \bbd_0\|^2}$ and decays exponentially.
Now we proceed by replacing the term $\E{\|\nabla f(\bbx_t) - \bbd_t\|^2}$ in \eqref{eq:proof_of_rate_thm_cvx_100} by its upper bound in \eqref{eq:done} and $\gamma_{t+1}$ by $2/(t+9)$ to write
{\small\begin{align}\label{eq:proof_of_rate_thm_cvx_200}
&\E{f(\bbx_{t+1}) -f(\bbx^*)}\leq\left(1-\frac{2}{t+9}\right)\E{(f(\bbx_{t}) -f(\bbx^*))}\nonumber\\&+\frac{R\sqrt{Q}}{(t+9)^{4/3}} +\frac{2LR^2}{(t+9)^{2}}.
\end{align}}
Note that we can write $(t+9)^2=(t+9)^{4/3} (t+9)^{2/3}\geq (t+9)^{4/3} 9^{2/3}\geq 4(t+9)^{4/3}$.Therefore, 
{\small\begin{align}\label{eq:proof_of_rate_thm_cvx_300}
&\E{f(\bbx_{t+1}) -f(\bbx^*)} \leq \left(1-\frac{2}{t+9}\right)\E{(f(\bbx_{t})-f(\bbx^*))}\nonumber\\&+\frac{2R\sqrt{Q}+LD^2/2}{(t+9)^{4/3}}.
\end{align}}
We use induction to prove for $t\geq 0$,
\begin{align*}
\E{f(\bbx_{t}) -f(\bbx^*)}\leq\frac{Q'}{(t+9)^{1/3}},
\end{align*}
where $Q'=\max\{9^{1/3}(f(\bbx_{0}) -f(\bbx^*)),2R\sqrt{2Q}+LR^2/2\}$. For $t=0$, we have that $\E{f(\bbx_{t}) -f(\bbx^*)}\leq\frac{Q'}{9^{1/3}}$, which is turn follows from the definition of $Q'$. Assume for the induction hypothesis holds for $t=k$. Then, for $t=k+1$, we have,
\begin{align*}
&\E{f(\bbx_{k+1}) -f(\bbx^*)} \leq \left(1-\frac{2}{k+9}\right)\E{(f(\bbx_{k})-f(\bbx^*))}\nonumber\\&+\frac{2R\sqrt{2Q}+LD^2/2}{(k+9)^{4/3}}\nonumber\\
&\le \left(1-\frac{2}{k+9}\right)\frac{Q'}{(t+9)^{1/3}}+\frac{Q'}{(t+9)^{4/3}} \le \frac{Q'}{(t+10)^{1/3}}.
\end{align*}
Thus, for $t\geq 0$ from Lemma \ref{lemma-estimation} we have that,
{\small\begin{align}\label{eq:claim:for_induction}
\E{f(\bbx_{t}) -f(\bbx^*)}\leq\frac{Q'}{(t+9)^{1/3}}=O\left(\frac{d^{1/3}}{(t+9)^{1/3}}\right).
\end{align}}
where $Q'=\max\{2(f(\bbx_{0}) -f(\bbx^*)),2R\sqrt{2Q}+LR^2/2\}$.
\end{proof}
\begin{proof}[Proof of Theorem \ref{th:2}(1)]
Then, we have,
{\small\begin{align}
	\label{eq:dual_rdsa}
	&F(\mathbf{x}_{t+1}) \leq F(\mathbf{x}_{t}) + \gamma_{t}\langle \mathbf{g}(\mathbf{x}_{t}), \mathbf{v}_{t}-\mathbf{x}_{t}\rangle\nonumber\\&+\gamma_{t}\langle\nabla F(\mathbf{x}_{t}) -\mathbf{g}(\mathbf{x}_{t}), \mathbf{v}_{t}-\mathbf{x}_{t}\rangle + \frac{LR^{2}\gamma_{t}^{2}}{2}\nonumber\\
&\Rightarrow F(\mathbf{x}_{t+1}) \leq F(\mathbf{x}_{t}) + \gamma_{t}\langle \mathbf{g}(\mathbf{x}_{t}), \argmin_{\mathbf{v}\in\mathcal{C}}\langle\mathbf{v},\nabla F(\mathbf{x}_{t})\rangle-\mathbf{x}_{t}\rangle\nonumber\\&+\gamma_{t}\langle\nabla F(\mathbf{x}_{t})-\mathbf{g}(\mathbf{x}_{t}), \mathbf{v}_{t}-\mathbf{x}_{t}\rangle + \frac{LR^{2}\gamma_{t}^{2}}{2}\nonumber\\
&\Rightarrow F(\mathbf{x}_{t+1}) \leq F(\mathbf{x}_{t}) + \gamma_{t}\langle \nabla F(\mathbf{x}_{t}), \argmin_{\mathbf{v}\in\mathcal{C}}\langle\mathbf{v},\nabla F(\mathbf{x}_{t})\rangle-\mathbf{x}_{t}\rangle\nonumber\\&+\gamma_{t}\langle\nabla F(\mathbf{x}_{t})-\mathbf{g}(\mathbf{x}_{t}), \mathbf{v}_{t}-\argmin_{\mathbf{v}\in\mathcal{C}}\langle\mathbf{v},\nabla F(\mathbf{x}_{t})\rangle\rangle + \frac{LR^{2}\gamma_{t}^{2}}{2}\nonumber\\
&\Rightarrow F(\mathbf{x}_{t+1}) \leq F(\mathbf{x}_{t})-\gamma_{t}\mathcal{G}\left(\mathbf{x}_{t}\right)\nonumber\\&+\gamma_{t}\langle\nabla F(\mathbf{x}_{t})-\mathbf{g}(\mathbf{x}_{t}), \mathbf{v}_{t}-\argmin_{\mathbf{v}\in\mathcal{C}}\langle\mathbf{v},\nabla F(\mathbf{x}_{t})\rangle\rangle + \frac{LR^{2}\gamma_{t}^{2}}{2}\nonumber\\
	&\Rightarrow\gamma_{t}\E{\mathcal{G}\left(\mathbf{x}_{t}\right)} \le \E{F(\mathbf{x}_{t})-F(\mathbf{x}_{t+1})}+\gamma_{t}R\frac{\sqrt{2Q}}{(t+8)^{1/3}}+ \frac{LR^{2}\gamma_{t}^{2}}{2}+\nonumber\\
	&\Rightarrow \E{\mathcal{G}\left(\mathbf{x}_{t}\right)} \le \E{\frac{t+7}{2}F(\mathbf{x}_{t})-\frac{t+8}{2}F(\mathbf{x}_{t+1})+\frac{1}{2}F(\mathbf{x}_{t})}+R\frac{\sqrt{2Q}}{(t+8)^{1/3}}+ \frac{LR^{2}\gamma_{t}}{2}\nonumber\\
	&\Rightarrow \sum_{t=0}^{T-1}\E{\mathcal{G}\left(\mathbf{x}_{t}\right)} \le \E{\frac{7}{2}F(\mathbf{x}_{0})-\frac{T+7}{2}F(\mathbf{x}_{T})+\sum_{t=0}^{T-1}\left(\frac{1}{2}F(\mathbf{x}_{t})}+R\frac{\sqrt{2Q}}{(t+8)^{1/3}}+ \frac{LR^{2}\gamma_{t}}{2}\right)\nonumber\\
	&\Rightarrow \sum_{t=0}^{T-1}\E{\mathcal{G}\left(\mathbf{x}_{t}\right)} \le \E{\frac{7}{2}F(\mathbf{x}_{0})-\frac{7}{2}F(\mathbf{x}^{\ast})}+\sum_{t=0}^{T-1}\left(\frac{1}{2}\left(F(\mathbf{x}_{t})-F(\mathbf{x}^{\ast})\right)+R\frac{\sqrt{2Q}}{(t+8)^{1/3}}+ \frac{LR^{2}\gamma_{t}}{2}\right)\nonumber\\
	&\Rightarrow \sum_{t=0}^{T-1}\E{\mathcal{G}\left(\mathbf{x}_{t}\right)} \le \frac{7}{2}F(\mathbf{x}_{0})-\frac{7}{2}F(\mathbf{x}^{\ast})+\sum_{t=0}^{T-1}\left(\frac{Q^{'}+R\sqrt{2Q}}{2(t+8)^{1/3}}+ \frac{LR^{2}}{(t+8)}\right)\nonumber\\
	&\Rightarrow T\E{\min_{t=0,\cdots,T-1}\mathcal{G}\left(\mathbf{x}_{t}\right)} \le \frac{7}{2}F(\mathbf{x}_{0})-\frac{7}{2}F(\mathbf{x}^{\ast})+LR^{2}\ln(T+7)+\frac{Q^{'}+R\sqrt{2Q}}{2}(T+7)^{2/3}\nonumber\\
	&\Rightarrow\E{\min_{t=0,\cdots,T-1}\mathcal{G}\left(\mathbf{x}_{t}\right)} \le \frac{7(F(\mathbf{x}_{0})-F(\mathbf{x}^{\ast}))}{2T}+\frac{LR^{2}\ln(T+7)}{T}+\frac{Q^{'}+R\sqrt{2Q}}{2T}(T+7)^{2/3}.
	\end{align}}
\end{proof}
\section{Proofs for Improvised RDSA}
\label{app:3}

\begin{proof}[Proof of Lemma \ref{le:1}(2)]
Following as in the proof of RDSA, we have,
{\small\begin{align}\label{proof:bound_on_grad_3001}
&\E{\|\nabla f(\bbx_t) - \bbd_t\|^2\mid\ccalF_t} \nonumber\\
&\le\rho_t^2\E{\|\nabla f(\bbx_t)-g(\bbx_t; \bby_t,\bbz_t)\|^2\mid\ccalF_t}  
\nonumber\\&+(1-\rho_t)^2\|\nabla f(\bbx_t)-\nabla f(\bbx_{t-1})\|^2\nonumber\\
&+(1-\rho_t)^2\|\nabla f(\bbx_{t-1}) - \bbd_{t-1} \|^2 \nonumber\\&+(1-\rho_t)^2\beta_t\|\nabla f(\bbx_{t-1}) - \bbd_{t-1}\|^2\nonumber\\&+\frac{(1-\rho_t)^2}{\beta_t}\|\nabla f(\bbx_t)-\nabla f(\bbx_{t-1}) \|^2\nonumber\\
&+2\rho_t (1-\rho_t)\frac{c_{t}^{2}}{m^{2}}\left\|L\mathbf{v}\left(\mathbf{x},c_{t}\right)\right\|^{2}+\rho_t (1-\rho_t)\left\|\nabla f(\bbx_t)-\nabla f(\bbx_{t-1})\right\|^{2}\nonumber\\
&+\rho_t (1-\rho_t)\left\|\nabla f(\bbx_{t-1}) - \bbd_{t-1}\right\|^{2}\nonumber\\
&\Rightarrow \E{\|\nabla f(\bbx_t) - \bbd_t\|^2}\le 2\rho_t^2\sigma^{2}+4\rho_{t}^{2}\E{\|\nabla F(\bbx_t,\bby_t)\|^{2}}\nonumber\\&+4\rho_{t}^{2}\E{\|g(\bbx_t; \bby_t,\bbz_t)\|^2} \nonumber\\
&+\left(1-\rho_t+\frac{(1-\rho_t)^2}{\beta_{t}}\right)\E{\|\nabla f(\bbx_t)-\nabla f(\bbx_{t-1})\|^2}
\nonumber\\&+\left(1-\rho_t+(1-\rho_t)^2\beta_{t}\right)\E{\|\nabla f(\bbx_{t-1}) - \bbd_{t-1} \|^2} \nonumber\\&
+\frac{\rho_t}{2} (1-\rho_t)c_{t}^{2}L^{2}M(\mu)\nonumber\\
&\le 2\rho_t^2\sigma^{2}+4\rho_{t}^{2}L_1^2+8\rho_{t}^{2}\left(1+\frac{s(d)}{m}\right)L_{1}^{2}+\left(\frac{1+m}{2m}\right)\rho_{t}^{2}c_{t}^{2}L^{2}M(\mu) \nonumber\\
&+\left(1-\rho_t+\frac{(1-\rho_t)^2}{\beta_{t}}\right)\E{\|\nabla f(\bbx_t)-\nabla f(\bbx_{t-1})\|^2}
\nonumber\\&+\left(1-\rho_t+(1-\rho_t)^2\beta_{t}\right)\E{\|\nabla f(\bbx_{t-1}) - \bbd_{t-1} \|^2} \nonumber\\&
+\frac{\rho_t}{2m^{2}}c_{t}^{2}L^{2}M(\mu),
\end{align}}
where we used the gradient approximation bounds as stated in \eqref{eq:grad_approx} and used Young's inequality to substitute the inner products and in particular substituted $2\langle \nabla f(\bbx_t)-\nabla f(\bbx_{t-1}) , \nabla f(\bbx_{t-1}) - \bbd_{t-1} \rangle$ by the upper bound $\beta_t \|\nabla f(\bbx_{t-1}) - \bbd_{t-1}\|^2+(1/\beta_t)\|\nabla f(\bbx_t)-\nabla f(\bbx_{t-1}) \|^2$ where $\beta_t>0$ is a free parameter. \\
According to Assumption \ref{as:4}, the norm $\|\nabla f(\bbx_t)-\nabla f(\bbx_{t-1}) \|$ is bounded above by $L\|\bbx_t-\bbx_{t-1}\|$. In addition, the condition in Assumption \ref{as:1} implies that $L\|\bbx_t-\bbx_{t-1}\|= L\gamma_{t}\|\bbv_t-\bbx_t\| \leq \gamma_{t}LR$. Therefore, we can replace $\|\nabla f(\bbx_t)-\nabla f(\bbx_{t-1}) \|$ by its upper bound $\gamma_{t}LR$ and since we assume that $\rho_t\leq1$ we can replace all the terms $(1-\rho_t)^2$. Furthermore, using $\beta_{t}:=\rho_t/2$ we have,
{\small\begin{align}\label{proof:bound_on_grad_5001}
&\E{\|\nabla f(\bbx_t) - \bbd_t\|^2} \nonumber\\
&\leq 2\rho_t^2\sigma^{2}+4\rho_t^2L_{1}^{2}+8\rho_t^2\left(1+\frac{s(d)}{m}\right)L_1^{2}+\frac{\rho_t}{2m^{2}}c_{t}^{2}L^{2}M(\mu)\nonumber\\
&+\gamma_{t}^2(1-\rho_t)\left(1+\frac{2}{\rho_t}\right)L^2R^2+\left(\frac{1+m}{2m}\right)\rho_t^2c_{t}^{2}L^{2}M(\mu)\nonumber\\
&+(1-\rho_t)\left(1+\frac{\rho_t}{2}\right)\E{\|\nabla f(\bbx_{t-1}) - \bbd_{t-1} \|^2} .
\end{align}}
Now using the inequalities $(1-\rho_t) (1+(2/\rho_t))\leq (2/\rho_t)$ and $(1-\rho_t)(1+({\rho_t}/{2}))\leq (1-\rho/2)$ we obtain
{\small\begin{align}\label{proof:bound_on_grad_7001}
&\E{\|\nabla f(\bbx_t) - \bbd_t\|^2}\leq 2\rho_t^2\sigma^{2}+4\rho_t^2L_{1}^{2}+8\rho_t^2\left(1+\frac{s(d)}{m}\right)L_1^{2}\nonumber\\
&+\left(\frac{1+m}{2m}\right)\rho_t^2c_{t}^{2}L^{2}M(\mu)+\frac{2L^2R^2\gamma_{t}^2 }{\rho_t}+\frac{\rho_t}{2m^{2}}c_{t}^{2}L^{2}M(\mu)\nonumber\\&+\left(1-\frac{\rho_t}{2}\right))\E{\|\nabla f(\bbx_{t-1}) - \bbd_{t-1} \|^2}.
\end{align}}
\end{proof}
\begin{proof}[Proof of Theorem \ref{th:1}(2)]
Now using the result in Lemma \ref{lemma-estimation} we can characterize the convergence of the sequence of expected errors $\E{\|\nabla f(\bbx_{t})-\bbd_t\|^2}$ to zero. To be more precise, using the result in Lemma \ref{le:1} and setting $\gamma_t=2/(t+8)$, $\rho_t=4/\left(1+\frac{d}{m}\right)^{1/3}(t+8)^{2/3}$ and $c_{t}=2\sqrt{m}/\sqrt{M(\mu)}(t+8)^{1/3}$, we have,
{\small\begin{align}\label{eq:almost_done1}
&\E{\|\nabla f(\bbx_t) - \bbd_t\|^2} 
\nonumber\\&\leq \left(1-\frac{2}{\left(1+\frac{d}{m}\right)^{1/3}(t+8)^{2/3}}\right)\E{\|\nabla F(\bbx_{t-1}) - \bbd_{t-1}\|^2}\nonumber\\&+\frac{32\left(1+\frac{d}{m}\right)^{-1/3}\sigma^2+64L_{1}^{2}\left(1+\frac{d}{m}\right)^{-1/3}+128\left(1+\frac{d}{m}\right)^{2/3}L_{1}^{2}}{(t+8)^{4/3}}\nonumber\\&+\frac{2L^2R^2\left(1+\frac{d}{m}\right)^{2/3}+416\left(1+\frac{d}{m}\right)^{2/3}L^{2}}{(t+8)^{4/3}}.
\end{align}}
According to the result in Lemma \ref{lemma-estimation}, the inequality in \eqref{eq:almost_done} implies that 
\begin{equation}\label{eq:done1}
 \E{\|\nabla f(\bbx_t) - \bbd_t\|^2} \leq \overline{Q}_{ir}+\frac{Q_{ir}}{(t+8)^{2/3}} \leq \frac{Q_{ir}}{(t+8)^{2/3}},
\end{equation}
where {\small$Q_{ir}=32\left(1+\frac{d}{m}\right)^{-1/3}\sigma^2+128\left(1+\frac{d}{m}\right)^{2/3}L_{1}^{2}+64\left(1+\frac{d}{m}\right)^{-1/3}L_{1}^{2}+2L^2R^2\left(1+\frac{d}{m}\right)^{2/3}+416\left(1+\frac{d}{m}\right)^{2/3}L^{2}$} and $\overline{Q}_{ir}$ is a function of $\E{\|\nabla f(\bbx_0) - \bbd_0\|^2}$ and decays exponentially.
Now we proceed by replacing the term $\E{\|\nabla f(\bbx_t) - \bbd_t\|^2}$ in \eqref{eq:proof_of_rate_thm_cvx_100} by its upper bound in \eqref{eq:done1} and $\gamma_{t+1}$ by $2/(t+9)$ to write
\begin{align}\label{eq:proof_of_rate_thm_cvx_2001}
&\E{f(\bbx_{t+1}) -f(\bbx^*)} 
\leq \left(1-\frac{2}{t+9}\right)\E{(f(\bbx_{t}) -f(\bbx^*))}\nonumber\\&+\frac{R\sqrt{2Q_{ir}}}{(t+9)^{4/3}} +\frac{2LR^2}{(t+9)^{2}}.
\end{align}
Note that we can write $(t+9)^2=(t+9)^{4/3} (t+9)^{2/3}\geq (t+9)^{4/3} 9^{2/3}\geq 4(t+9)^{4/3}$.Therefore, 
\begin{align}\label{eq:proof_of_rate_thm_cvx_3001}
&\E{f(\bbx_{t+1}) -f(\bbx^*)} 
\leq \left(1-\frac{2}{t+9}\right)\E{(f(\bbx_{t})-f(\bbx^*))}\nonumber\\&+ \frac{2R\sqrt{Q}+LD^2/2}{(t+9)^{4/3}} .
\end{align}
Following the induction steps as in \eqref{eq:claim:for_induction}, we have,
\begin{align}\label{eq:claim:for_induction1}
\E{f(\bbx_{t}) -f(\bbx^*) } \leq \frac{Q_{ir}^{'}}{(t+8)^{1/3}} = O\left(\frac{(d/m)^{1/3}}{(t+9)^{1/3}}\right).
\end{align}
where $Q_{ir}'=\max\{2(f(\bbx_{0}) -f(\bbx^*)),2R\sqrt{2Q_{ir}}+LR^2/2\}$.
\end{proof}
\begin{proof}[Proof of Theorem \ref{th:2}(2)]
Following as in \eqref{eq:dual_rdsa}, we have,
{\small\begin{align}
	\label{eq:dual_irdsa}
	&\gamma_{t}\E{\mathcal{G}\left(\mathbf{x}_{t}\right)} \le \E{F(\mathbf{x}_{t})-F(\mathbf{x}_{t+1})}+\gamma_{t}R\frac{\sqrt{2Q_{ir}}}{(t+8)^{1/3}}+ \frac{LR^{2}\gamma_{t}^{2}}{2}+\nonumber\\
	&\Rightarrow \E{\mathcal{G}\left(\mathbf{x}_{t}\right)} \le \E{\frac{t+7}{2}F(\mathbf{x}_{t})-\frac{t+8}{2}F(\mathbf{x}_{t+1})+\frac{1}{2}F(\mathbf{x}_{t})}+R\frac{\sqrt{2Q_{ir}}}{(t+8)^{1/3}}+ \frac{LR^{2}\gamma_{t}}{2}\nonumber\\
	&\Rightarrow \sum_{t=0}^{T-1}\E{\mathcal{G}\left(\mathbf{x}_{t}\right)} \le \E{\frac{7}{2}F(\mathbf{x}_{0})-\frac{T+7}{2}F(\mathbf{x}_{T})+\sum_{t=0}^{T-1}\left(\frac{1}{2}F(\mathbf{x}_{t})}+R\frac{\sqrt{2Q_{ir}}}{(t+8)^{1/3}}+ \frac{LR^{2}\gamma_{t}}{2}\right)\nonumber\\
	&\Rightarrow \sum_{t=0}^{T-1}\E{\mathcal{G}\left(\mathbf{x}_{t}\right)} \le \E{\frac{7}{2}F(\mathbf{x}_{0})-\frac{7}{2}F(\mathbf{x}^{\ast})}+\sum_{t=0}^{T-1}\left(\frac{1}{2}\left(F(\mathbf{x}_{t})-F(\mathbf{x}^{\ast})\right)+R\frac{\sqrt{2Q_{ir}}}{(t+8)^{1/3}}+ \frac{LR^{2}\gamma_{t}}{2}\right)\nonumber\\
	&\Rightarrow \sum_{t=0}^{T-1}\E{\mathcal{G}\left(\mathbf{x}_{t}\right)} \le \frac{7}{2}F(\mathbf{x}_{0})-\frac{7}{2}F(\mathbf{x}^{\ast})+\sum_{t=0}^{T-1}\left(\frac{Q_{ir}^{'}+R\sqrt{2Q_{ir}}}{2(t+8)^{1/3}}+ \frac{LR^{2}}{(t+8)}\right)\nonumber\\
	&\Rightarrow T\E{\min_{t=0,\cdots,T-1}\mathcal{G}\left(\mathbf{x}_{t}\right)} \le \frac{7}{2}F(\mathbf{x}_{0})-\frac{7}{2}F(\mathbf{x}^{\ast})+LR^{2}\ln(T+7)+\frac{Q_{ir}^{'}+R\sqrt{2Q_{ir}}}{2}(T+7)^{2/3}\nonumber\\
	&\Rightarrow\E{\min_{t=0,\cdots,T-1}\mathcal{G}\left(\mathbf{x}_{t}\right)} \le \frac{7(F(\mathbf{x}_{0})-F(\mathbf{x}^{\ast}))}{2T}+\frac{LR^{2}\ln(T+7)}{T}+\frac{Q_{ir}^{'}+R\sqrt{2Q_{ir}}}{2T}(T+7)^{2/3}
	\end{align}}
\end{proof}
\section{Proofs for KWSA}
\label{app:4}
\begin{proof}[Proof of Lemma \ref{le:1}(3)]
Following as in the proof of Lemma \ref{le:1}, we have,
{\small\begin{align}\label{proof:bound_on_grad_3002}
&\E{\|\nabla f(\bbx_t) - \bbd_t\|^2}\nonumber\\
&\le(1-\rho_t)^2\E{\|\nabla f(\bbx_t)-\nabla f(\bbx_{t-1})\|^2} 
\nonumber\\&+(1-\rho_t)^2\|\E{\nabla f(\bbx_{t-1}) - \bbd_{t-1} \|^2} \nonumber\\&+(1-\rho_t)^2\beta_t\E{\|\nabla f(\bbx_{t-1}) - \bbd_{t-1}\|^2}\nonumber\\&+\frac{(1-\rho_t)^2}{\beta_t}\E{\|\nabla f(\bbx_t)-\nabla f(\bbx_{t-1}) \|^2}\nonumber\\
&+\frac{\rho_t}{2}(1-\rho_t)c_{t}^{2}L^{2}d\nonumber\\&+\rho_t(1-\rho_t)\E{\left\|\nabla f(\bbx_t)-\nabla f(\bbx_{t-1})\right\|^{2}}\nonumber\\
&+\rho_t (1-\rho_t)\E{\left\|\nabla f(\bbx_{t-1}) - \bbd_{t-1}\right\|^{2}}\nonumber\\
&\le 2\rho_t^2\sigma^{2}+2\rho_{t}^{2}c_{t}^{2}dL^{2} \nonumber\\
&+\left(1-\rho_t+\frac{(1-\rho_t)^2}{\beta_{t}}\right)\E{\|\nabla f(\bbx_t)-\nabla f(\bbx_{t-1})\|^2}
\nonumber\\&+\left(1-\rho_t+(1-\rho_t)^2\beta_{t}\right)\E{\|\nabla f(\bbx_{t-1}) - \bbd_{t-1} \|^2} \nonumber\\&
+\frac{\rho_t}{2} (1-\rho_t)c_{t}^{2}L^{2}d\nonumber\\
&\le 2\rho_t^2\sigma^{2}+2\rho_{t}c_{t}^{2}dL^2
\nonumber\\&+\left(1-\rho_t+\frac{(1-\rho_t)^2}{\beta_{t}}\right)\E{\|\nabla f(\bbx_t)-\nabla f(\bbx_{t-1})\|^2}\nonumber\\
&+\left(1-\rho_t+(1-\rho_t)^2\beta_{t}\right)\E{\|\nabla f(\bbx_{t-1}) - \bbd_{t-1} \|^2},
\end{align}}
where we used the gradient approximation bounds as stated in \eqref{eq:grad_approx} and used Young's inequality to substitute the inner products and in particular substituted $2\langle \nabla f(\bbx_t)-\nabla f(\bbx_{t-1}) , \nabla f(\bbx_{t-1}) - \bbd_{t-1} \rangle$ by the upper bound $\beta_t \|\nabla f(\bbx_{t-1}) - \bbd_{t-1}\|^2+(1/\beta_t)\|\nabla f(\bbx_t)-\nabla f(\bbx_{t-1}) \|^2$ where $\beta_t>0$ is a free parameter. \\
According to Assumption \ref{as:4}, the norm $\|\nabla f(\bbx_t)-\nabla f(\bbx_{t-1}) \|$ is bounded above by $L\|\bbx_t-\bbx_{t-1}\|$. In addition, the condition in Assumption \ref{as:1} implies that $L\|\bbx_t-\bbx_{t-1}\|= L\gamma_{t}\|\bbv_t-\bbx_t\| \leq \gamma_{t}LR$. Therefore, we can replace $\|\nabla f(\bbx_t)-\nabla f(\bbx_{t-1}) \|$ by its upper bound $\gamma_{t}LR$ and since we assume that $\rho_t\leq1$ we can replace all the terms $(1-\rho_t)^2$. Furthermore, using $\beta_{t}:=\rho_t/2$ we have,
\begin{align}\label{proof:bound_on_grad_5002}
	&\E{\|\nabla f(\bbx_t) - \bbd_t\|^2} \nonumber\\
	&\leq 2\rho_t^2\sigma^{2}+2\rho_{t}c_{t}^{2}dL^2+\gamma_{t}^2(1-\rho_t)\left(1+\frac{2}{\rho_t}\right)L^2R^2\nonumber\\
	&+(1-\rho_t)\left(1+\frac{\rho_t}{2}\right)\E{\|\nabla f(\bbx_{t-1}) - \bbd_{t-1} \|^2} .
\end{align}
	Now using the inequalities $(1-\rho_t) (1+(2/\rho_t))\leq (2/\rho_t)$ and $(1-\rho_t)(1+({\rho_t}/{2}))\leq (1-\rho/2)$ we obtain
\begin{align}\label{proof:bound_on_grad_7002}
	&\E{\|\nabla f(\bbx_t) - \bbd_t\|^2}\leq 2\rho_t^2\sigma^{2}+2\rho_{t}c_{t}^{2}dL^2+\frac{2L^2R^2\gamma_{t}^2 }{\rho_t}\nonumber\\&+\left(1-\frac{\rho_t}{2}\right))\E{\|\nabla f(\bbx_{t-1}) - \bbd_{t-1} \|^2}.
\end{align}
\end{proof}
\begin{proof}[Proof of Theorem \ref{th:1}(3)]
Now using the result in Lemma \ref{lemma-estimation} we can characterize the convergence of the sequence of expected errors $\E{\|\nabla f(\bbx_{t})-\bbd_t\|^2}$ to zero. To be more precise, using the result in Lemma \ref{le:1} and setting $\gamma_t=2/(t+8)$, $\rho_t=4/(t+8)^{2/3}$ and $c_{t}=2/\sqrt{d}(t+8)^{1/3}$ for any $\epsilon > 0$ to obtain
\begin{align}\label{eq:almost_done2}
& \E{\|\nabla f(\bbx_t) - \bbd_t\|^2} \leq \nonumber\\&\left(1-\frac{2}{(t+8)^{2/3}}\right)\E{\|\nabla F(\bbx_{t-1}) - \bbd_{t-1}\|^2}\nonumber\\&+\frac{32\sigma^2+32L^{2}+2L^2R^2}{(t+8)^{4/3}}.
\end{align}
According to the result in Lemma \ref{lemma-estimation}, the inequality in  \eqref{eq:almost_done} implies that 
\begin{equation}\label{eq:done2}
\E{\|\nabla f(\bbx_t) - \bbd_t\|^2} \leq \frac{Q_{kw}}{(t+8)^{2/3}},
\end{equation}
where 
{\small\begin{align*}
Q=\max\left\{4\|\nabla f(\bbx_0) - \bbd_0\|^2,32\sigma^2+32L^{2}+2L^2R^2\right\}
\end{align*}}
Now we proceed by replacing the term $\E{\|\nabla f(\bbx_t) - \bbd_t\|^2}$ in \eqref{eq:proof_of_rate_thm_cvx_100} by its upper bound in \eqref{eq:done1} and $\gamma_{t+1}$ by $2/(t+9)$ to write
\begin{align}\label{eq:proof_of_rate_thm_cvx_2002}
&\E{f(\bbx_{t+1}) -f(\bbx^*)} 
\leq \left(1-\frac{2}{t+9}\right)\E{(f(\bbx_{t}) -f(\bbx^*))}\nonumber\\&+ \frac{R\sqrt{Q_{kw}}}{(t+9)^{4/3}}+\frac{2LR^2}{(t+9)^{2}}.
\end{align}
Note that we can write $(t+9)^2=(t+9)^{4/3} (t+9)^{2/3}\geq (t+9)^{4/3} 9^{2/3}\geq 4(t+9)^{4/3}$.Therefore, 
\begin{align}\label{eq:proof_of_rate_thm_cvx_3002}
&\E{f(\bbx_{t+1}) -f(\bbx^*)} \leq \left(1-\frac{2}{t+9}\right)\E{(f(\bbx_{t})-f(\bbx^*))}\nonumber\\&+\frac{2R\sqrt{Q_{kw}}+LD^2/2}{(t+9)^{4/3}} .
\end{align}
Thus, for $t\geq 0$ by induction we have,
\begin{align}\label{eq:claim:for_induction2}
\E{f(\bbx_{t}) -f(\bbx^*) } \leq \frac{Q'}{(t+9)^{1/3}} = O\left(\frac{d^{0}}{(t+9)^{1/3}}\right).
\end{align}
where $Q'=\max\{2(f(\bbx_{0}) -f(\bbx^*)),2R\sqrt{Q_{kw}}+LR^2/2\}$.
\end{proof}
\begin{proof}[Proof of Theorem \ref{th:2}(3)]
Following as in \eqref{eq:dual_rdsa}, we have,
{\small\begin{align}
	\label{eq:dual_kwsa}
	&\gamma_{t}\E{\mathcal{G}\left(\mathbf{x}_{t}\right)} \le \E{F(\mathbf{x}_{t})-F(\mathbf{x}_{t+1})}+\gamma_{t}R\frac{\sqrt{2Q_{kw}}}{(t+8)^{1/3}}+ \frac{LR^{2}\gamma_{t}^{2}}{2}+\nonumber\\
	&\Rightarrow \E{\mathcal{G}\left(\mathbf{x}_{t}\right)} \le \E{\frac{t+7}{2}F(\mathbf{x}_{t})-\frac{t+8}{2}F(\mathbf{x}_{t+1})+\frac{1}{2}F(\mathbf{x}_{t})}+R\frac{\sqrt{2Q_{kw}}}{(t+8)^{1/3}}+ \frac{LR^{2}\gamma_{t}}{2}\nonumber\\
	&\Rightarrow \sum_{t=0}^{T-1}\E{\mathcal{G}\left(\mathbf{x}_{t}\right)} \le \E{\frac{7}{2}F(\mathbf{x}_{0})-\frac{T+7}{2}F(\mathbf{x}_{T})+\sum_{t=0}^{T-1}\left(\frac{1}{2}F(\mathbf{x}_{t})}+R\frac{\sqrt{2Q_{kw}}}{(t+8)^{1/3}}+ \frac{LR^{2}\gamma_{t}}{2}\right)\nonumber\\
	&\Rightarrow \sum_{t=0}^{T-1}\E{\mathcal{G}\left(\mathbf{x}_{t}\right)} \le \E{\frac{7}{2}F(\mathbf{x}_{0})-\frac{7}{2}F(\mathbf{x}^{\ast})}+\sum_{t=0}^{T-1}\left(\frac{1}{2}\left(F(\mathbf{x}_{t})-F(\mathbf{x}^{\ast})\right)+R\frac{\sqrt{2Q_{kw}}}{(t+8)^{1/3}}+ \frac{LR^{2}\gamma_{t}}{2}\right)\nonumber\\
	&\Rightarrow \sum_{t=0}^{T-1}\E{\mathcal{G}\left(\mathbf{x}_{t}\right)} \le \frac{7}{2}F(\mathbf{x}_{0})-\frac{7}{2}F(\mathbf{x}^{\ast})+\sum_{t=0}^{T-1}\left(\frac{Q_{kw}^{'}+R\sqrt{2Q_{kw}}}{2(t+8)^{1/3}}+ \frac{LR^{2}}{(t+8)}\right)\nonumber\\
	&\Rightarrow T\E{\min_{t=0,\cdots,T-1}\mathcal{G}\left(\mathbf{x}_{t}\right)} \le \frac{7}{2}F(\mathbf{x}_{0})-\frac{7}{2}F(\mathbf{x}^{\ast})+LR^{2}\ln(T+7)+\frac{Q_{kw}^{'}+R\sqrt{2Q_{kw}}}{2}(T+7)^{2/3}\nonumber\\
	&\Rightarrow\E{\min_{t=0,\cdots,T-1}\mathcal{G}\left(\mathbf{x}_{t}\right)} \le \frac{7(F(\mathbf{x}_{0})-F(\mathbf{x}^{\ast}))}{2T}+\frac{LR^{2}\ln(T+7)}{T}+\frac{Q_{kw}^{'}+R\sqrt{2Q_{kw}}}{2T}(T+7)^{2/3}
	\end{align}}
\end{proof}
\section{Proofs for Non Convex Stochastic Frank Wolfe}
\label{app:5}
\begin{proof}[Proof of Theorem \ref{th:4}]
We reuse the following characterization derived earlier:
\begin{Lemma}
	\label{le:2nc}
	Let Assumptions \ref{as:3}-\ref{as:6} hold. Given the recursion in \eqref{eq:2}, we have that $\|\nabla f(\bbx_t) - \bbd_t\|^2$ satisfies
	\begin{align}
	\label{eq:le2nc}
	&\E{\|\nabla f(\bbx_t) - \bbd_t\|^2}\leq 2\rho_t^2\sigma^{2}+4\rho_t^2L_{1}^{2}\nonumber\\&+8\rho_t^2\left(1+\frac{s(d)}{m}\right)L_1^{2}+\left(\frac{1+m}{2m}\right)\rho_t^2c_{t}^{2}L^{2}M(\mu)\nonumber\\
	&+\frac{2L^2R^2\gamma^2}{\rho_t}+\frac{\rho_t}{2m^{2}}c_{t}^{2}L^{2}M(\mu)\nonumber\\&+\left(1-\frac{\rho_t}{2}\right))\E{\|\nabla f(\bbx_{t-1}) - \bbd_{t-1} \|^2}.
	\end{align}
\end{Lemma}
Now using the result in Lemma \ref{lemma-estimation} we can characterize the convergence of the sequence of expected errors $\E{\|\nabla f(\bbx_{t})-\bbd_t\|^2}$ to zero. To be more precise, using the result in Lemma \ref{le:1} and setting $\gamma=T^{-3/4}$, $\rho_t=4/\left(1+\frac{d}{m}\right)^{1/3}(t+8)^{1/2}$ and $c_{t}=2\sqrt{m}/\sqrt{M(\mu)}(t+8)^{1/4}$ to obtain for all $t=0,\cdots,T-1$,
{\small\begin{align}\label{eq:almost_done4}
 & \E{\|\nabla f(\bbx_t) - \bbd_t\|^2} \nonumber\\
 &\leq \left(1-\frac{2}{\left(1+\frac{d}{m}\right)^{1/3}(t+8)^{1/2}}\right)\E{\|\nabla F(\bbx_{t-1}) - \bbd_{t-1}\|^2}\nonumber\\&+\frac{32\sigma^2+64L_{1}^{2}+128\left(1+\frac{d}{m}\right)^{1/3}L_{1}^{2}}{(t+8)}\nonumber\\&+\frac{8L^2R^2\left(1+\frac{d}{m}\right)^{1/3}+416L^{2} }{(t+8)}.
 \end{align}}
Using Lemma \ref{lemma-estimation}, we then have,
{\small\begin{align}
 \label{eq:nc_bound}
 \E{\|\nabla f(\bbx_t) - \bbd_t\|^2} = O\left(\frac{(d/m)^{2/3}}{(t+9)^{1/2}}\right), \forall~t=0,\cdots,T-1
 \end{align}}
Finally, we have,
{\small\begin{align}
 \label{eq:dual_nc}
 &F(\mathbf{x}_{t+1}) \leq F(\mathbf{x}_{t}) + \gamma_{t}\langle \mathbf{d}_{t}, \mathbf{v}_{t}-\mathbf{x}_{t}\rangle\nonumber\\&+\gamma\langle\nabla F(\mathbf{x}_{t}) -\mathbf{d}_{t}, \mathbf{v}_{t}-\mathbf{x}_{t}\rangle + \frac{LR^{2}\gamma^{2}}{2}\nonumber\\
 &\leq F(\mathbf{x}_{t}) + \gamma\langle \mathbf{d}_{t}, \argmin_{\mathbf{v}\in\mathcal{C}}\langle\mathbf{v},\nabla F(\mathbf{x}_{t})\rangle-\mathbf{x}_{t}\rangle\nonumber\\&+\gamma\langle\nabla F(\mathbf{x}_{t})-\mathbf{d}_{t}, \mathbf{v}_{t}-\mathbf{x}_{t}\rangle + \frac{LR^{2}\gamma^{2}}{2}\nonumber\\
 &\leq F(\mathbf{x}_{t}) + \gamma\langle \nabla F(\mathbf{x}_{t}), \argmin_{\mathbf{v}\in\mathcal{C}}\langle\mathbf{v},\nabla F(\mathbf{x}_{t})\rangle-\mathbf{x}_{t}\rangle\nonumber\\&+\gamma\langle\nabla F(\mathbf{x}_{t})-\mathbf{d}_{t}, \mathbf{v}_{t}-\argmin_{\mathbf{v}\in\mathcal{C}}\langle\mathbf{v},\nabla F(\mathbf{x}_{t})\rangle\rangle + \frac{LR^{2}\gamma^{2}}{2}\nonumber\\
 &\leq F(\mathbf{x}_{t})-\gamma\mathcal{G}\left(\mathbf{x}_{t}\right) + \frac{LR^{2}\gamma^{2}}{2}\nonumber\\
 &+\gamma\langle\nabla F(\mathbf{x}_{t})-\mathbf{d}_{t}, \mathbf{v}_{t}-\argmin_{\mathbf{v}\in\mathcal{C}}\langle\mathbf{v},\nabla F(\mathbf{x}_{t})\rangle\rangle\nonumber\\
 &\Rightarrow \gamma\E{\mathcal{G}\left(\mathbf{x}_{t}\right)} \le \E{F(\mathbf{x}_{t})}-\E{F(\mathbf{x}_{t+1})}\nonumber\\&+\gamma R\E{\left\|\nabla F(\mathbf{x}_{t})-\mathbf{d}_{t}\right\|} +\frac{LR^{2}\gamma^{2}}{2}\nonumber\\
 &\le \E{F(\mathbf{x}_{t})}-\E{F(\mathbf{x}_{t+1})}+\gamma_{t}R\sqrt{\E{\left\|\nabla F(\mathbf{x}_{t})-\mathbf{d}_{t}\right\|^{2}}} +\frac{LR^{2}\gamma^{2}}{2}\nonumber\\
 &\le \E{F(\mathbf{x}_{t})}-\E{F(\mathbf{x}_{t+1})}+Q_{nc}\gamma\rho_{t}^{1/2}R(d/m)^{1/3} +\frac{LR^{2}\gamma^{2}}{2}\nonumber\\
 &\Rightarrow \E{\mathcal{G}_{min}}T\gamma \le \E{F(\mathbf{x}_{0})}-\E{F(\mathbf{x}_{t+1})}\nonumber\\&+Q_{nc}\gamma R(d/m)^{1/3}\sum_{t=0}^{T-1}\rho_{t}^{1/2} +\frac{LR^{2}T\gamma^{2}}{2}\nonumber\\
 &\Rightarrow \E{\mathcal{G}_{min}} \le \frac{\E{F(\mathbf{x}_{0})}-\E{F(\mathbf{x}^{*})}}{T\gamma}\nonumber\\&+\gamma Q_{nc}R(d/m)^{1/3}\frac{\sum_{t=0}^{T-1}\rho_{t}^{1/2}}{T\gamma} +\frac{LR^{2}T\gamma^{2}}{2T\gamma}\nonumber\\
 &\Rightarrow \E{\mathcal{G}_{min}} \le \frac{\E{F(\mathbf{x}_{0})}-\E{F(\mathbf{x}^{\ast})}}{T^{1/4}}\nonumber\\&+\frac{Q_{nc}Rd^{1/3}}{T^{1/4}m^{1/3}} +\frac{LR^{2}}{2T},
\end{align}}
 where $\mathcal{G}_{min}=\min_{t=0,\cdots,T-1}\mathcal{G}\left(\mathbf{x}_{t}\right)$.
 \end{proof}
\end{document}